\documentclass[twoside, 12pt]{article}
\usepackage{amssymb,mathrsfs,psfrag,graphicx}
\usepackage{CJK,bbm,amsfonts,mathrsfs,amsthm,amssymb,amsmath,indentfirst,color}
\usepackage{hyperref}
\def\v{\varepsilon}

\def\t{\theta}
\def\T{\Theta}
\def\k{\kappa}

\def\a{\alpha}
\def\b{\beta}
\def\g{\gamma}
\def\d{\delta}
\def\l{\lambda}
\def\f{\frac}
\def\p{\phi}

\def\ra{\rightarrow}

\def\z{\zeta}

\def\di{\displaystyle}
\def\i{\infty}
 \newtheorem{lemma}{\bf Lemma}[section]
       \newtheorem{theorem}[lemma]{\bf Theorem}

       \newtheorem{proposition}[lemma]{\bf Proposition}
       \newtheorem{remark}[lemma]{\bf Remark}

\def\Xint#1{\mathchoice
{\XXint\displaystyle\textstyle{#1}}%
{\XXint\textstyle\scriptstyle{#1}}%
{\XXint\scriptstyle\scriptscriptstyle{#1}}%
{\XXint\scriptscriptstyle\scriptscriptstyle{#1}}%
\!\int}
\def\XXint#1#2#3{{\setbox0=\hbox{$#1{#2#3}{\int}$}
\vcenter{\hbox{$#2#3$}}\kern-.5\wd0}}

\def\dashint{\Xint-}

\def\Xparallel#1{\mathchoice
{\XXparallel\displaystyle\textstyle{#1}}%
{\XXparallel\textstyle\scriptstyle{#1}}%
{\XXparallel\scriptstyle\scriptscriptstyle{#1}}%
{\XXparallel\scriptscriptstyle\scriptscriptstyle{#1}}%
\!\parallel}
\def\XXparallel#1#2#3{{\setbox0=\hbox{$#1{#2#3}{\parallel}$}
\vcenter{\hbox{$#2#3$}}\kern-.6\wd0}}
\def\dpr{\Xparallel-}

\begin{document}
\date{}
\title{\Large \bf Zero dissipation limit of full compressible Navier-Stokes equations with Riemann initial data}
\author{\small \textbf{Feimin Huang},\thanks{F. Huang is supported was supported in part by NSFC
Grant No. 10825102 for distinguished youth scholar, and National
Basic Research Program of China (973 Program) under Grant No.
2011CB808002. E-mail: fhuang@amt.ac.cn.}\quad
  \textbf{Song Jiang}\thanks{S. Jiang is supported by NSFC Grant No. 40890154
and the National Basic Research Program under the Grant
2011CB309705. E-mail: jiang@iapcm.ac.cn.}\quad and \textbf{Yi
Wang}\thanks{Corresponding author. Y. Wang is supported by NSFC grant No. 10801128 and No.
11171326. E-mail: wangyi@amss.ac.cn.}} \maketitle \small
 $^\ast$ $^\ddag$ Institute of Applied Mathematics, AMSS, and
Hua Loo-Keng Key Laboratory of Mathematics, CAS, Beijing 100190,
China

 $^{\dag}$ Institute of Applied Physics and Computational
Mathematics, Beijing 100088, China\\

 {\bf Abstract:} We consider the zero dissipation limit of the full compressible
Navier-Stokes equations with Riemann initial data in the case of
superposition of two rarefaction waves and a contact discontinuity.
It is proved that for any suitably small viscosity $\v$ and heat
conductivity $\k$ satisfying the relation \eqref{viscosity}, there
exists a unique global piecewise smooth solution to the compressible
Navier-Stokes equations. Moreover, as the viscosity $\v$ tends to
zero, the Navier-Stokes solution converges uniformly to the Riemann
solution of superposition of two rarefaction waves and a contact
discontinuity to the corresponding Euler equations with the same
Riemann initial data away from the initial line $t=0$ and the
contact discontinuity located at $x=0$.

\renewcommand{\theequation}{\arabic{section}.\arabic{equation}}
\setcounter{equation}{0}
\section{Introduction}
 We study the zero dissipation limit of the solution to the Navier-Stokes
 equations of a compressible heat-conducting gas in Lagrangian coordinate:
\begin{equation}
 \left\{
  \begin{array}{l}
    \di v_t-u_x=0,\\
    \di u_t+p_x=\v(\f{u_x}{v})_x,\\
    \di (e+\f{u^2}{2})_t+(pu)_x=\k(\f{\t_x}{v})_x+\v(\f{uu_x}{v})_x
  \end{array}
\right. \label{NS}
\end{equation}
with Riemann initial data
\begin{equation}
(v,u,\t)(0,x)=\left\{
\begin{array}{ll}
(v_-,u_-,\t_-),&x<0,\\
(v_+,u_+,\t_+),&x>0,
\end{array}
\right.\label{DD}
\end{equation}
 where the functions
$v(x,t)>0,u(x,t),\t(x,t)>0$ represent the specific volume, velocity
and the absolute temperature of the gas, respectively. And
$p=p(v,\t)$ is the pressure, $e=e(v,\t)$ is the internal energy,
$\v>0$ is the viscosity constant and $\k>0$ is the coefficient of
heat conduction. Here we consider an ideal and polytropic gas, that is
\begin{equation*}
p=\f{R\t}{v},\qquad e=\f{R\t}{\g-1}, \label{state-equa}
\end{equation*}
with $\g>1, R>0$ being gas constants.

The study of the asymptotic behavior of viscous flows, as the
viscosity tends to zero, is one of the important problems in the
theory of compressible fluid flows. When the solution of the
inviscid flow is smooth, the zero dissipation limit problem can be
solved by classical scaling method. However, the inviscid
compressible flow contains discontinuities, such as shock waves, in
general. In this case, it is also conjectured that a general weak entropy solution
to the inviscid flow should be the strong limit of the solution to
the corresponding viscous flows with the same initial data as the
viscosity vanishes.

It is well known that the solution to the Riemann problem for
the Euler equations consists of three basic wave patterns, that is,
shock, rarefaction wave and contact discontinuity. Moreover, the
Riemann solution is essential in the theory for the Euler equations
as it captures both local and global behavior of general solutions.

For hyperbolic conservation laws with the uniform viscosity
$$
u_t+f(u)_x=\v u_{xx},
$$
where $f(u)$ satisfies some assumptions to ensure the hyperbolic
nature of the corresponding inviscid system, Goodman-Xin \cite{GX}
verified the limit for piecewise smooth solutions separated by
non-interacting shock waves using a matched asymptotic expansion
method. Later, Yu \cite{Y} proved it for hyperbolic conservation
laws with both shock and initial layers. In 2005, important progress
made by Bianchini-Bressan\cite{BB} justifies the vanishing
viscosity limit in BV-space even though the problem is still
unsolved for the physical system such as the compressible Navier-Stokes
equations.

For the compressible isentropic Navier-Stokes equations
 where the conservation of energy in
\eqref{NS} is neglected in the isentropic regime, Hoff-Liu
\cite{HL} firstly proved the vanishing viscosity limit for a
piecewise constant shock with initial layer. Later, Xin
\cite{X} justified the limit for rarefaction waves. Then,
Wang \cite{WH} generalized the result of Goodman-Xin
\cite{GX} to the isentropic Navier-Stokes equations.

Recently, Chen-Perepelitsa \cite{CP} proved the convergence of the isentropic
compressible Navier-Stokes equations to the compressible Euler equations as the viscosity
vanishes in Eulerian coordinates for general initial data by using compensated
compactness method if the far field does not contain vacuum.
Note that this result allows the initial data containing vacuum in the
interior domain. However, the framework of compensated compactness
is basically limited to $2\times 2$ systems so far, so that this
result could not be applied to the full compressible Navier-Stokes
equations (\ref{NS}).

For the full compressible Navier-Stokes equations, there are investigations on the
 limits to the Euler system for the basic wave patterns in the literature.
We refer to Jiang-Ni-Sun \cite{JNS} and Xin-Zeng \cite{XZ} for the
rarefaction wave, Wang \cite{W} for the shock wave, Ma \cite{M} for
the contact discontinuity and Huang-Wang-Yang \cite{HWY, HWY1} for
the superposition of two rarefaction waves and a contact
discontinuity and the superposition of rarefaction and shock waves,
respectively. We should point out that the limit shown in
\cite{JNS} was  for the discontinuous initial data while the other results mentioned were  for (well-prepared) smooth data.

In this paper, we shall investigate the zero dissipation limit of the full
Navier-Stokes equations \eqref{NS} with Riemann initial data
\eqref{DD} in the case of the superposition of two rarefaction waves
and a contact discontinuity. The local and global
well-posedness of the full system \eqref{NS} or the corresponding
isentropic system with discontinuous initial data is systematically
studied by Hoff, etc., see \cite{H1, H2, H3, H4, H5, H6, CHT}.
In order to get the zero dissipation limit to the Riemann
solution of the Euler system, we shall combine the local existence of
solutions with discontinuous data from \cite{H3} and the
time-asymptotic stability analysis to the compressible Navier-Stokes
equations \eqref{NS1}. Compared with the previous result
\cite{HWY} where the same limit process is studied for
(well-prepared) smooth initial data, the main difficulty in the proof here lies
in the discontinuity of the initial data.
The discontinuity of the initial data for the volume $v(t,x)$ will propagate for all the time along
the particle path due to the hyperbolic regime while the smoothing effects will
also be performed on the velocity $u(t,x)$ and the temperature $\t(t,x)$ by the parabolic structure,
and this interaction of the discontinuity and smoothing effects brings technical difficulties.
To circumvent such difficulties, we shall choose suitable weight functions to carry out the weighted
energy estimates in terms of the superposition wave structure (see Remark 3.7), and use the
energy method of Huang-Li-Matsumura \cite{HLM} for the stability of
two rarefaction waves with a contact discontinuity in the middle,
where the authors obtained a new estimate on the heat kernel which
can be applied to the study of the stability of the viscous contact wave in the framework
of the rarefaction wave (see Lemma \ref{lemma3}). Namely, the anti-derivative
variable of the perturbation is not necessary and the estimates to
the perturbation itself are also available to get the stability of
the viscous contact wave.

Without loss of generality, we assume the following relation between the
viscosity constant $\v$ and the heat-conducing coefficient $\k$ of
system \eqref{NS} as in \cite{JNS}:
\begin{equation}
\left\{
\begin{array}{l}
\di \k=O(\v),\qquad \qquad \rm as \qquad\v\rightarrow0;\\
\di \nu\doteq\f{\k(\v)}{\v}\geq c>0\qquad {\rm for ~some ~positive~
constant}~ c,\quad \rm as \quad\v\rightarrow0.
\end{array}
\right. \label{viscosity}
\end{equation}

If $\k=\v=0$ in \eqref{NS}, then the corresponding Euler system reads as
\begin{equation}
\left\{
\begin{array}{l}
v_t-u_x=0,\\
u_t+p_x=0,\\[2mm]
\displaystyle{\Big( e+\f{u^2}{2}\Big)_t+(pu)_x=0.}
\end{array}
\right.\label{euler}
\end{equation}
It can be easily computed that the eigenvalues of the Jacobi
matrix of the flux function to \eqref{euler} are
\begin{equation}
\l_1=-\sqrt{\f{\g p}{v}},\quad \l_2=0,\quad \l_3=\sqrt{\f{\g p}{v}}.
\label{eigen}
\end{equation}
It is well known that the first and third
characteristic fields of \eqref{euler} are genuinely nonlinear and the second one
is linearly degenerate (see\cite{S}).

For the Euler equations, we know that there are three basic
wave patterns, shock, rarefaction wave and contact discontinuity.
And the Riemann solution to the Euler equations has a basic wave
pattern consisting the superposition of these three waves with the
contact discontinuity in the middle. For later use, let us firstly
recall the wave curves for the two types of basic waves studied in
this paper.

Given the right end state $(v_+,u_+,\t_+)$ with $v_+, \t_+>0$,  the following wave
curves in the phase space $\{(v,u,\t)|v>0, \t>0\}$ are defined for the Euler equations.

$\bullet$ Contact discontinuity curve:
\begin{equation}
CD(v_+,u_+,\t_+)= \{(v,u,\t)  |  u=u_+, p=p_+, v \not\equiv v_+
 \}. \label{CD}
\end{equation}

$\bullet$ $i$-Rarefaction wave curve $(i=1,3)$:
\begin{equation}
 R_i (v_+, u_+, \theta_+):=\Bigg{ \{} (v, u, \theta)\Bigg{ |}u<u_+ ,~u=u_+-\int^v_{v_+}
 \lambda_i(\eta,
s_+) \,d\eta,~ s(v, \theta)=s_+\Bigg{ \}},\label{Ri}
\end{equation}
where $s_+=s(v_+,\t_+)$ and $\l_i=\l_i(v,s)$ defined in \eqref{eigen} is the $i$-th
characteristic speed of the Euler system \eqref{euler}.

Now, we define the solution profile that consists of the superposition of two rarefaction waves
and a contact discontinuity. Let $ (v_-,u_-,\t_-) \in$
$R_1$-$CD$-$R_3(v_+,u_+,\t_+)$. Then, there exist uniquely two
intermediate states $(v_*, u_*,\t_*)$ and $(v^*, u^*,\t^*)$, such
that $(v_*, u_*,\t_*)\in R_1(v_-, u_-,\t_-)$, $(v_*, u_*,\t_*)\in
CD(v^*, u^*,\t^*)$  and $(v^*, u^*,\t^*)\in R_3(v_+,u_+,\t_+)$.

Thus, the  wave pattern $(\bar V,\bar U,\bar\T)(t,x)$ consisting of
1-rarefaction wave, 2-contact discontinuity and 3-rarefaction wave
that solves the corresponding
 Riemann problem of the Euler system \eqref{euler} can be defined by
\begin{eqnarray}
 \left(\begin{array}{cc} \bar V\\ \bar U \\
 \bar  \T
\end{array}
\right)(t, x)= \left(\begin{array}{cc}v^{r_1}+ v^{cd}+ v^{r_3}\\ u^{r_1}+ u^{cd}+ u^{r_3} \\
\t^{r_1}+ \t^{cd}+ \t^{r_3}
\end{array}
\right)(t, x) -\left(\begin{array}{cc} v_*+v^*\\ u_*+u^*\\
\t_*+\t^*
\end{array}
\right) ,\label{RS}
\end{eqnarray}
where $(v^{r_1}, u^{r_1}, \t^{r_1} )(t,x)$ is the 1-rarefaction wave
defined in \eqref{Ri} with the right state $(v_+, u_+, \t_+)$
replaced by $(v_*, u_*, \theta_* )$, $(v^{cd}, u^{cd}, \t^{cd}
)(t,x)$ is the contact discontinuity   defined in \eqref{CD} with
the states $(v_-, u_-, \t_-)$ and $(v_+, u_+, \t_+)$ replaced by
$(v_*, u_*, \theta_* )$ and $(v^*, u^*, \theta^* )$ respectively,
and $(v^{r_3}, u^{r_3}, \t^{r_3})(t,x)$ is the 3-rarefaction wave
defined in \eqref{Ri} with the left state $(v_-, u_-, \t_-)$
replaced by $(v^*, u^*, \theta^* )$.

Now we state the main result as follows.
\begin{theorem}\label{limit-th} Given a Riemann solution $(\bar V,\bar
U,\bar\T)(t,x)$ defined in \eqref{RS}, which is superposition of
two rarefaction waves and a contact discontinuity for the Euler
system \eqref{euler}, there exist small positive constants
$\delta_0$ and $\v_0$, such that if $\varepsilon\leq\varepsilon_0$ and the wave strength
$\delta\doteq|(v_+-v_-,u_+-u_-,\t_+-\t_-)|\leq \delta_0$,
 then the compressible Navier-Stokes equations \eqref{NS} with \eqref{DD} and
 \eqref{viscosity} admits a unique global piece-wise smooth solution
$(v^\v,u^\v,\t^\v)(t,x)$ satisfying that
\begin{itemize}
\item The quantities $u^\v,\t^\v$, $p(v^\v,\t^\v)-\v\f{u^\v_x}{v^\v}$ and $\f{\t^\v_x}{v^\v}$ are
continuous for $t>0$, and the jumps in $v^\v,u^\v_x,\t^\v_x$ at
$x=0$ satisfies
\begin{equation*}
|([v^\v(t,0)],[u^\v_x(t,0)],[\t^\v_x(t,0)])|\leq Ce^{-\f{ct}{\v}},
\end{equation*}
where the constants $C$ and $c$ are independent of $t$ and $\v$.

\item Moreover, under the condition \eqref{viscosity}, it holds that
\begin{equation}
\lim_{\v\rightarrow
0}\sup_{(t,x)\in\Sigma_h}|(v^\v,u^\v,\t^\v)(t,x)-(\bar V,\bar
U,\bar\T)(t,x)|=0,\quad \forall h>0, \label{limit}
\end{equation}
where $\Sigma_h=\big\{(t,x)|t\geq h, \f{x}{\sqrt{\v+t}}\geq h
\v^{\a},0\leq\a<\f12\big\}$.
\end{itemize}
\end{theorem}
\begin{remark} Theorem \ref{limit-th} shows that, away from the
initial time $t=0$ and the contact discontinuity located at $x=0$, there exists a unique
global solution $(v^\v,u^\v,\t^\v)(t,x)$  of the compressible
Navier-Stokes equations \eqref{NS} which converges to the Riemann
solution $(\bar V,\bar U,\bar\T)(t,x)$ consisting of two rarefaction
waves and a contact discontinuity when $\v$ and $\k$ satisfy the relation \eqref{viscosity}
and $\v$ tends to zero. Moreover, the convergence is uniform on the
set $\Sigma_h$ for any $h>0$.
\end{remark}
\noindent{\bf Notations.} In the paper, we always use the notation
$\di \dashint_{\mathbf R}=\int_{\mathbf{R}^+}+\int_{\mathbf{R}^-}$,
$\|\cdot\|$ to denote the usual $L^2(\mathbf{R})$ norm, $\|\cdot
\dpr$ to denote the piecewise $L^2$ norm, that is, $\di \|f
\dpr^2=\dashint_{\mathbf R}f^2dy$. $\|\cdot\|_1$ and
$\|\cdot\dpr_1$ represent the $H^1(\mathbf{R})$ norm and piece-wise
$H^1(\mathbf{R}^\pm)$ norm, respectively. And the notation $[\cdot]$
represents the jump of the function $\cdot$ at $x=0$ or $y=0$ if without confusion.

\section{Approximate profiles}
\setcounter{equation}{0}
Introduce the following scaled variables
\begin{equation}
y=\f{x}{\v},\quad \tau=\f{t}{\v}, \label{scaling}
\end{equation}
and set
\begin{equation*}\label{new-unknown}
(v^\v,u^\v,\t^\v)(t,x)=(v,u,\t)(\tau,y).
\end{equation*}
Then the new unknown functions $(v,u,\t)(\tau,y)$ satisfies the system
\begin{equation}\label{NS1}
 \left\{
  \begin{array}{l}
    \di v_\tau-u_y=0,\\
    \di u_\tau+p_y=(\f{u_y}{v})_y,\\
    \di
     \f{R}{\g-1}\t_\tau+pu_y=\nu(\f{\t_y}{v})_y+\f{u^2_y}{v},
  \end{array}
\right.
\end{equation}
with the scaled heat conductivity $\nu=\f{\k}{\v}$ in \eqref{viscosity} satisfying
$$
\nu_0\leq \nu\leq \nu_1, ~~{\rm uniformly}~{\rm in}~\v ~{\rm as}~ \v\rightarrow0+,~{\rm for}~ {\rm some}
~ {\rm positive}~ {\rm constants}~ \nu_0~ {\rm and}~ \nu_1. $$

Note that the Riemann solution $(\bar V,\bar U,\bar\T)(t,x)$ in \eqref{RS} is invariant under the scaling
transformation \eqref{scaling}, thus to prove the limit \eqref{limit} in Theorem \ref{limit-th},
it is sufficient to show the following limit
\begin{equation}
\lim_{\v\rightarrow
0}\sup_{(\tau,y)\in\Sigma^1_h}|(v,u,\t)(\tau,y)-(\bar V,\bar
U,\bar\T)(\tau,y)|=0,\quad \forall h>0, \label{limit-1}
\end{equation}
where $\Sigma_h^1$ is the corresponding region of $\Sigma_h$ in the new coordinates $(\tau,y)$ defined by
\begin{equation}\label{sigma-1}
\Sigma^1_h=\Big\{(\tau,y)|\tau\geq \f h\v, \f{y}{\sqrt{1+\tau}}\geq \f {h}{\v^{\f12-\a}}, 0\leq \a<\f12\Big\}.
\end{equation}

Now we study the Navier-Stokes equations \eqref{NS1}. The corresponding wave profiles to \eqref{CD} and
\eqref{Ri} can be defined approximately as follows. We start from the viscous contact wave to \eqref{CD}.

\subsection{Viscous contact wave}

If $(v_-,u_-,\t_-)\in CD(v_+,u_+,\t_+)$, i.e.,
$$
u_-=u_+,~p_-=p_+,~v_-\neq v_+,
$$
then the Riemann problem, that is, the Euler system \eqref{euler}
with Riemann initial data
$$
(v,u,\t)(\tau=0,y)=\left\{
\begin{array}{ll}
(v_-,u_-,\t_-),\qquad & y<0,\\
(v_+,u_+,\t_+),\qquad & y>0,
\end{array}
\right.
$$
admits a single contact discontinuity solution
\begin{equation*}
(v^{cd},u^{cd},\t^{cd})(\tau,y)=\left\{
\begin{array}{ll}
(v_-,u_+,\t_-),\qquad & y<u_+\tau,~ \tau>0,\\
(v_+,u_+,\t_+),\qquad & y>u_+\tau,~ \tau>0.
\end{array}
\right.
\end{equation*}

As in \cite{HMX}, the viscous version of the above
contact discontinuity, called viscous contact wave
$(V^{CD},U^{CD},\T^{CD})(\tau,y)$, can be defined as follows. Since it is
expected that
$$
P^{CD}\approx p_+=p_-, \quad {\rm and}\quad |U^{CD}-u_+|\ll1,
$$
the leading order of the energy equation $\eqref{NS1}_3$ is
$$
\f{R}{\g-1}\T_\tau+p_+U_y=\nu(\f{\T_y}{V})_y.
$$
Then, similar to \cite{HLM} or \cite{HWY}, one can get the following nonlinear diffusion equation
$$
\T_\tau=a \Big(\f{\T_y}{\T}\Big)_y,\quad \T(\tau,\pm)=\t_\pm,\quad a=\f{\nu
p_+(\g-1)}{R^2\g}.
$$
The above diffusion equation has a unique self-similar solution
$\hat\T(\tau,y)=\hat\T(\f{y}{\sqrt{1+\tau}})$.

Thus, the viscous contact wave $(V^{CD},U^{CD},\T^{CD})(\tau,y)$ can be defined by
\begin{equation}
\begin{array}{ll}
\di V^{CD}(\tau,y)=\f{R\hat\T(\tau,y)}{p_+},\\[4mm]
\di U^{CD}(\tau,y)=u_+
+\f{\nu(\g-1)}{R\g}\f{\hat\T_{y}(\tau,y)}{\hat\T(\tau,y)},\\[5mm]
\di \T^{CD}(\tau,y)=\hat\T(\tau,y)+\f{R\g-\nu(\g-1)}{\g p_+}\hat\T_\tau.
\end{array}
\label{Viscous-CD}
\end{equation}
Here, it is straightforward to check that the viscous contact wave
defined in \eqref{Viscous-CD} satisfies
\begin{equation}
|\hat\T-\t_\pm|+(1+\tau)^{\f12}|\hat\T_y|+(1+\tau)|\hat\T_{yy}|
=O(1)\d^{CD} e^{-\f{c_0y^2}{1+\tau}}, \quad\mbox{as }|y|\ra+\i , \label{CD-P}
\end{equation}
where $\d^{CD}=|\t_+-\t_-|$ represents the strength
of the viscous contact wave and $c_0$ is a positive constant. Note that in \eqref{Viscous-CD},
the higher order term $\f{R\g-\nu(\g-1)}{\g p_+}\hat\T_\tau$ is introduced in
$\Theta^{CD}(\tau,y)$ to make the viscous contact wave
$(V^{CD},U^{CD},\T^{CD})(\tau,y)$ satisfy the momentum equation
exactly. Correspondingly,  $(V^{CD},U^{CD},\T^{CD})(\tau,y)$ satisfies the system
\begin{equation}
\left\{
\begin{array}{l}
\di V^{\scriptscriptstyle CD}_{\tau}-U^{CD}_{y}=0,\\[2mm]
\di U^{CD}_{\tau}+P^{CD}_{y}=\Big(\f{U^{CD}_{y}}{V^{CD}}\Big)_y,\\[4mm]
\di
\f{R}{\g-1}\T^{CD}_{\tau}+P^{CD}U^{CD}_{y}=\nu\Big(\f{\T^{CD}_{y}}{V^{CD}}\Big)_y+\f{(U^{CD}_{y})^2}{V^{CD}}+Q^{CD},
\end{array}
\right. \label{CD-system}
\end{equation}
where $\di P^{CD}=\f{R\T^{CD}}{V^{CD}}$ and the error term $Q^{CD}$
satisfies
\begin{equation}
\di  Q^{CD} =O(1)\d^{CD}(1+\tau)^{-2}e^{-\f{c_0y^2}{1+\tau}},\qquad {\rm as}~~|y|\ra +\i,
\label{Q-CD}
\end{equation}
for some positive constant $c_0$.

\subsection{Approximate rarefaction waves}
We now turn to the approximate rarefaction waves to \eqref{Ri}. Since there is no exact
rarefaction wave profile for the Navier-Stokes equations, the
following approximate rarefaction wave profile, which satisfies the Euler
equations, is motivated by \cite{X}. For the completeness of
presentation, we include its definition and the properties in this
subsection.

If $(v_-, u_-, \theta_-) \in R_i (v_+, u_+, \theta_+), (i=1,3)$, then
there exists an $i$-rarefaction wave $(v^{r_i}, u^{r_i},
\t^{r_i})(y/\tau)$ which is a global  solution of the following Riemann problem:
\begin{eqnarray}
\left\{
\begin{array}{l}
\di  v_{\tau}- u_{y}= 0,\\[1mm]
\di u_{\tau} +  p_{y}(v, \theta) = 0 , \\[2mm]
\di \frac{R}{\g-1}\theta_\tau + p(v, \theta) u_y =0,\\[1mm]
\di  (v, u, \t)(0,y)=\left\{
\begin{array}{l}
\di (v_-, u_-, \t_-), ~~~~  y<
0 ,\\
\di  (v_+, u_+, \t_+),~~~~ y> 0 .
\end{array}
\right.
\end{array} \right.\label{rarefaction}
\end{eqnarray}
Consider the following inviscid Burgers equation with Riemann data:
\begin{equation}
\left\{
\begin{array}{l}
w_\tau+ww_y=0,\\[2mm]
w(\tau=0,y)=\left\{
\begin{array}{ll}
w_-,&y<0,\\
w_+,&y>0.
\end{array}
 \right.
  \end{array}
 \right.\label{Burgers}
\end{equation}
If $w_-<w_+$, then the Riemann problem \eqref{Burgers} admits a rarefaction
wave solution
\begin{equation}
w^r(\tau,y)=w^r(\f y\tau)=\left\{
\begin{array}{ll}
w_-,&\f y\tau\leq w_-,\\[1mm]
\f y\tau,&w_-\leq \f y\tau\leq w_+,\\[1mm]
w_+,&\f y\tau\geq w_+.
\end{array}
\right.\label{B-R}
\end{equation}
Thus, the Riemann solution in \eqref{rarefaction} can be expressed explicitly
through the above rarefaction wave \eqref{B-R} to the Burgers equation, that is,
\begin{eqnarray}
\left\{
\begin{array}{l}
\di  s^{r_i}(\tau,y)=s(v^{r_i}(\tau,y),\t^{r_i}(\tau,y))=s_+, \\[2mm]
\di w_\pm=\l_{i\pm}:=\l_i(v_\pm,\t_\pm),  \\[1mm]
\di w^r(\f y\tau)= \l_i(v^{r_i}(\tau,y),s_+), \\[1mm]
\di u^{r_i}(\tau,y)=u_+-\int^{v^{r_i}(\tau,y)}_{v_+}  \l_i(v,s_+) dv.
\end{array} \right.\label{AR-NS-1}
\end{eqnarray}

In order to construct the approximate rarefaction wave $(V^{R_i}, U^{R_i},
\Theta^{R_i})(\tau,y)$ corresponding to \eqref{Ri}, we first consider the
following approximate rarefaction wave to the Burgers equation:
\begin{eqnarray}
\left\{
\begin{array}{l}
\di w_{\tau}+ww_{y}=0,\\
\di w( 0,y )=w_0(y)=\f{w_++w_-}{2}+\f{w_+-w_-}{2}\tanh y.
\end{array}
\right.\label{AB}
\end{eqnarray}
Note that the solution $w^R(\tau,y)$ of the problem \eqref{AB} is
given by
$$
w^R(\tau,y)=w_0(x_0(\tau,y)),\qquad x=x_0(\tau,y)+w_0(x_0(\tau,y))\tau.
$$
And $w^R(\tau,y)$ has the following properties, the proof of which can be
found in \cite{MN1,X}:
 \vskip 2mm
\begin{lemma}\label{lemma-R}
Let $w_-<w_+$, then $\eqref{AB}$ has a unique smooth solution $w^R(\tau,y)$
satisfying
\begin{enumerate}
\item[(1)] $w_-< w^R(\tau,y)<w_+,~(w^R)_y(\tau,y)> 0 $;

\item[(2)] For any $1\leq p\leq +\infty$, there exists a constant $C$ such that
$$
\begin{array}{ll}
\| \f{\partial}{\partial y}w^R(\tau,\cdot)\|_{L^p(\mathbf{R})}\leq
C\min\big{\{}(w_+-w_-),~
(w_+-w_-)^{1/p}\tau^{-1+1/p}\big{\}}, \\[2mm]
 \| \f{\partial^2}{\partial y^2}w^R(\tau,\cdot)\|_{L^p(\mathbf{R})}\leq
C\min\big{\{}(w_+-w_-),~ \tau^{-1}\big{\}};
\end{array}
$$
\item[(3)] If $ y-w_-\tau<0$, then
$$
\begin{array}{l}
 |w^R(\tau,y)-w_-|\leq (w_+-w_-)e^{-2|y-w_-\tau|},\\[2mm]
 |\f{\partial}{\partial y}w^R(\tau,y)|\leq
2(w_+-w_-)e^{-2|y-w_-\tau|};
\end{array}
$$
 If $ y-w_+\tau> 0$, then
$$
\begin{array}{l}
 |w^R(\tau,y)-w_+|\leq (w_+-w_-)e^{-2|y-w_+\tau|},\\[2mm]
 |\f{\partial}{\partial x}w^R(\tau,y)|\leq
2(w_+-w_-)e^{-2|y-w_+\tau|};
\end{array}
$$

\item[(4)]$\sup\limits_{y\in\mathbf{R}} |w^R(\tau,y)-w^r(\f y\tau)|\leq
\min\big\{w_+-w_-,\f1{\tau}\ln(1+\tau)\big\}$.
\end{enumerate}
\end{lemma}

Then, corresponding to \eqref{AR-NS-1}, the approximate rarefaction wave profile denoted by
$(V^{R_i}, U^{R_i}, \T^{R_i}) (\tau,y)~(i=1,3)$ to \eqref{Ri} can be defined  by
\begin{eqnarray}
\left\{
\begin{array}{l}
\di  S^{R_i}(\tau,y)=s(V^{R_i}(\tau,y),\T^{R_i}(\tau,y))=s_+,\\[1mm]
\di w_\pm=\l_{i\pm}:=\l_i(v_\pm,\t_\pm), \\[2mm]
\di w^R(1+\tau,y)= \l_i(V^{R_i}(\tau,y),s_+),\\[1mm]
\di U^{R_i}(\tau,y)=u_+-\int^{V^{R_i}(\tau,y)}_{v_+}  \l_i(v,s_+) dv.
\end{array} \right.\label{AR-NS}
\end{eqnarray}
Note that $(V^{R_i}, U^{R_i}, \Theta^{R_i})(\tau,y)$ defined above
satisfies
\begin{equation}
\left\{
\begin{array}{ll}
\di  V^{R_i}_\tau-U^{R_i} _{y} = 0, \\[1mm]
 \di    U^{R_i}_\tau+P^{R_i}_y
    =0,\\[2mm]
 \di \f{R}{\g-1}  \T^{R_i}_\tau
        + P^{R_i} U^{R_i}_y
    =0,
\end{array}
\right.\label{R-system}
\end{equation}
where $ P^{R_i}=p( V^{R_i}, \T  ^{R_i})$.

By virtue of Lemmas \ref{lemma-R}, the properties on the approximate
rarefaction waves $(V^{R_i}, U^{R_i}, \Theta^{R_i})(\tau,y)$ can be summarized as follows.
\begin{lemma}\label{lemma-R-NS}
 The approximate rarefaction waves
$(V^{R_i}, U^{R_i}, \Theta^{R_i})(\tau,y)~(i=1,3)$ constructed in
\eqref{AR-NS} have the following properties:
\begin{enumerate}
\item[(1)] $U^{R_i}_x(\tau,y)>0$ for $y\in \mathbf{R}$, $\tau>0$;
\item[(2)] For any $1\leq p\leq +\i,$ the following estimates holds,
$$
\begin{array}{ll}
\|(V^{R_i},U^{R_i}, \Theta^{R_i})_y\|_{L^p(dy)} \leq
C\min\big{\{}\d^{R_i},~
(\d^{R_i})^{1/p}(1+\tau)^{-1+1/p}\big{\}},\\[2mm]
\|(V^{R_i},U^{R_i}, \Theta^{R_i})_{yy}\|_{L^p(dy)} \leq
C\min\big{\{}\d^{R_i},~ (1+\tau)^{-1}\big{\}},\\
\end{array}
$$
where $\d^{R_i}=|(v_+,v_-, u_+,u_-, \t_+,\t_-)|$ is the $i$-rarefaction wave strength
and the positive constant $C$ is independent of $\tau$, but may only depend on $p$ and the wave strength;
\item[(3)] If $y\geq \l_{1+}(1+\tau)$, then
$$
\begin{array}{l}
 |(V^{R_1},U^{R_1},\T^{R_1})(\tau,y)-(v_-,u_-,\t_-)|\leq C\d^{R_1}e^{-2|y-\l_{1+}(1+\tau)|},\\[2mm]
 |(V^{R_1},U^{R_1},\T^{R_1})_y(\tau,y)|\leq
C \d^{R_1}e^{-2|y-\l_{1+}(1+\tau)|};
\end{array}
$$
If $y\leq \l_{3-}(1+\tau)$, then
$$
\begin{array}{l}
 |(V^{R_3},U^{R_3},\T^{R_3})(\tau,y)-(v_+,u_+,\t_+)|\leq C\d^{R_3}e^{-2|y-\l_{3-}(1+\tau)|},\\[2mm]
 |(V^{R_3},U^{R_3},\T^{R_3})_y(\tau,y)|\leq
C\d^{R_3} e^{-2|y-\l_{3-}(1+\tau)|};
\end{array}
$$

\item[(4)] There exists a positive constant $C$, such that
for all $\tau>0,$
$$
\sup_{y\in\mathbf{R}}|(V^{R_i},U^{R_i},
\Theta^{R_i})(\tau,y)-(v^{r_i},u^{r_i}, \theta^{r_i})(\f y\tau)|\leq
\f{C}{1+\tau}\ln(1+\tau).
$$
\end{enumerate}
\end{lemma}

\subsection{Superposition of rarefaction waves and contact
discontinuity}

Corresponding to \eqref{RS}, the approximate wave pattern $(V,U,\T)(\tau,y)$ of
the compressible Navier-Stokes equations \eqref{NS1} can be defined by
\begin{eqnarray}
 \left(\begin{array}{cc} V\\ U \\
  \T
\end{array}
\right)(\tau,y)= \left(\begin{array}{cc} V^{R_1}+ V^{CD}+ V^{R_3}\\ U^{R_1}+ U^{CD}+ U^{R_3} \\
\T^{R_1}+ \T^{CD}+ \T^{R_3}
\end{array}
\right)(\tau,y) -\left(\begin{array}{cc} v_*+v^*\\ u_*+u^*\\
\t_*+\t^*
\end{array}
\right) ,\label{sup-wave}
\end{eqnarray}
where $(V^{R_1}, U^{R_1}, \Theta^{R_1} )(\tau,y)$ is the approximate
1-rarefaction wave defined in \eqref{AR-NS} with the right state
$(v_+, u_+, \t_+)$ replaced by $(v_*, u_*, \theta_* )$, $(V^{CD},
U^{CD}, \Theta^{CD} )(\tau,y)$ is the viscous contact wave   defined in
\eqref{Viscous-CD} with the states $(v_-, u_-, \t_-)$ and $(v_+, u_+,
\t_+)$ replaced by $(v_*, u_*, \theta_* )$ and $(v^*, u^*, \theta^*
)$ respectively, and $(V^{R_3}, U^{R_3}, \Theta^{R_3} )(\tau,y)$ is the
approximate 3-rarefaction wave defined in \eqref{AR-NS} with the left
state $(v_-, u_-, \t_-)$ replaced by $(v^*, u^*, \theta^* )$.

Thus, from the properties of the viscous contact wave in \eqref{CD-P} and the approximate
rarefaction wave in Lemma 2.3, we
have the following relation between the approximate wave pattern
$(V,U,\T)(\tau,y)$ and the exact inviscid wave pattern
 $(\bar V,\bar U,\bar\T)(\tau,y)$ of the Euler equations
\begin{equation}
\di |(V,U,\T)(\tau,y)-(\bar V,\bar U,\bar\T)(\tau,y)| \di
\leq \f{C}{1+\tau}\ln(1+\tau)+C\d^{CD}e^{-\f{cy^2}{1+\tau}}.  \label{profile-s}
\end{equation}
Hence, to prove the zero dissipation limit \eqref{limit-1} on the set $\Sigma_h^1$ defined
in \eqref{sigma-1}, it is sufficient to show the following time-asymptotic behavior of the
solution to \eqref{NS1} around the approximate wave profile \eqref{sup-wave}, i.e.,
\begin{equation}\label{con}
\lim_{\tau\rightarrow+\i}\|(v,u,\t)(\tau,\cdot)-(V,U,\T)(\tau,\cdot)\|_{L^\i}=0.
\end{equation}

First, by \eqref{CD-system} and \eqref{R-system}, the superposition wave profile
$(V,U,\T)(\tau,y)$ defined in \eqref{sup-wave} satisfies the following system
\begin{equation*}
\left\{
\begin{array}{ll}
 \di V_\tau-U _{y} = 0, \\
  \di  U_\tau+P_y
    = (\frac{U_{y}}{V}) _{y}+Q_1,\\
 \di    \f{R}{\g-1} \T _\tau+P U_y
    =\nu( \frac{\T_{y}}{ V})_y+  \frac{U_y^2}{ V}
     +Q_2,
\end{array}
\right.
\end{equation*}
 where $ P =p( V , \T )$ and
$$\begin{array}{ll}
   \di Q_1&\di=(P-P^{R_1}-P^{CD}-P^{R_3})_y-\left(\f{U_y}{V}-\f{U^{CD}_y}{V^{CD}}\right)_y,\\
  \di Q_2&\di= (PU_y-P^{R_1}U^{R_1}_y-P^{CD}U^{CD}_y-P^{R_3}U^{R_3}_y)
  -\nu\left(\f{\T_y}{V}-\f{\T^{CD}_y}{V^{CD}}\right)_y\\
  &\di-\left(\f{U_y^2}{ V}-\f{(U^{CD}_y)^2}{ V^{CD}}\right)-Q^{CD}.
\end{array}$$
A direct calculation shows that
\begin{equation}
\begin{array}{lll}
\di Q_1&=&\di O(1)\Big\{|(V^{R_1}_y,\T^{R_1}_y)||(V^{CD}-v_*,\T^{CD}-\t_*,V^{R_3}-v^*,\T^{R_3}-\t^*)|\\[2mm]
&&\di
+|(V^{R_3}_y,\T^{R_3}_y)||(V^{R_1}-v_*,\T^{R_1}-\t_*,V^{CD}-v^*,\T^{CD}-\t^*)|\\[2mm]
&&\di+|(V^{CD}_y,\T^{CD}_y,U^{CD}_{y})||(V^{R_1}-v_*,\T^{R_1}-\t_*,V^{R_3}-v^*,\T^{R_3}-\t^*)|\\[2mm]
&& \di +|(U^{CD}_y,V^{CD}_y)||(U^{R_1}_y,V^{R_1}_y,U^{R_3}_y,V^{R_3}_y)|
+|(U^{R_1}_y,V^{R_1}_y)||(U^{R_3}_y,V^{R_3}_y)|\Big\}\\[2mm]
&&\di+O(1)\Big\{|U^{R_1}_{yy}|+|U^{R_3}_{yy}|+|U^{R_1}_y||V^{R_1}_y|+|U^{R_3}_y||V^{R_3}_y|\Big\}\\[2mm]
& :=&\di Q_{11}+Q_{12}.
\end{array}
\label{Q1}
\end{equation}
Similarly, we have
\begin{equation}
\begin{array}{lll}
\di Q_2&=&\di O(1)\Big\{|U^{R_1}_y||(V^{CD}-v_*,\T^{CD}-\t_*,V^{R_3}-v^*,\T^{R_3}-\t^*)|\\[2mm]
&&\di +|U^{R_3}_y||(V^{R_1}-v_*,\T^{R_1}-\t_*,V^{CD}-v^*,\T^{CD}-\t^*)|\\[2mm]
&&\di+|(U^{CD}_{y},V^{CD}_y,\T^{CD}_y)||(V^{R_1}-v_*,\T^{R_1}-\t_*,V^{R_3}-v^*,\T^{R_3}-\t^*)|\\[2mm]
&&\di +|(U^{CD}_y,V^{CD}_y,\T^{CD}_y)||(U^{R_1}_y,V^{R_1}_y,\T^{R_1}_y,U^{R_3}_y,V^{R_3}_y,\T^{R_1}_y)|
\\[2mm]
&&\di +|(U^{R_1}_y,V^{R_1}_y,\T^{R_1}_y)||(U^{R_3}_y,V^{R_3}_y,\T^{R_3}_y)|\Big\}\\[2mm]
&&\di +O(1)\Big\{|\T^{R_1}_{yy}|+|\T^{R_3}_{yy}|+|(U^{R_1}_y,V^{R_1}_y,
\T^{R_1}_y,U^{R_3}_y,V^{R_3}_y,\T^{R_3}_y)|^2\Big\}+|Q^{CD}|\\[2mm]
& :=&\di Q_{21}+Q_{22}+|Q^{CD}|.
\end{array}
\label{Q2}
\end{equation}
Here $Q_{11}$ and $Q_{21}$ represent the wave interaction terms coming from
the wave patterns in the different family, $Q_{12}$ and $Q_{22}$ stand for the error
terms due to the inviscid approximate  rarefaction wave profiles, and
$Q^{CD}$ is the error term defined in \eqref{Q-CD} due to the viscous
contact wave.

In fact, one can estimate the interaction terms $Q_{11}$ and $Q_{21}$ by
dividing the whole domain
$\Omega=\{(\tau,y)|(\tau,y)\in\mathbf{R}\times\mathbf{R}\}$ into three
regions:
\begin{eqnarray*} &&
\Omega_-=\{(\tau,y)\; |\; 2y\leq \l_{1*}(1+\tau)\}, \\
&& \Omega_{CD}=\{(\tau,y)\; |\; \l_{1*}(1+\tau)<2y<\l_3^*(1+\tau)\}, \\
&& \Omega_+=\{(\tau,y)\; |\; 2y\geq \l_3^*(1+\tau)\},
\end{eqnarray*}
where $\l_{1*}=\l_1(v_*,\t_*)$ and $\l_3^*=\l_3(v^*,\t^*)$.
Then, in each section the following estimates follow from
\eqref{CD-P} and Lemma \ref{lemma-R-NS}.
\begin{itemize}
\item In $\Omega_-$,
\begin{eqnarray*} &&
|(V^{R_3}-v^*, V^{R_3}_y)| = O(1)\d^{R_3}e^{-2\{|y|+|\l_3^*|(1+\tau)\}}, \\
&& |(V^{CD}-v_*,V^{CD}-v^*,V^{CD}_y)| = O(1)\d^{CD}e^{-\f{C\{|\l_{1*}|(1+\tau)\}^2}{1+\tau}}
=O(1)\d^{CD}e^{-C(1+\tau)}; \end{eqnarray*}

\item In $\Omega_{CD}$,
\begin{eqnarray*} &&
|(V^{R_1}-v_*, V^{R_1}_y) | =O(1)\d^{R_1} e^{-2\{|y|+|\l_{1*}|(1+\tau)\}}, \\
&& |(V^{R_3}-v^*, V^{R_3}_y)|=O(1)\d^{R_3} e^{-2\{|x|+|\l_3^*|(1+\tau)\}};
\end{eqnarray*}
\item In $\Omega_+$,
\begin{eqnarray*}  && |(V^{R_1}-v_*, V^{R_1}_y)|=O(1)\d^{R_1} e^{-2\{|x|+|\l_{1*}|(1+\tau)\}}, \\
&& |(V^{CD}-v_*, V^{CD}-v^*, V^{CD}_y)|=O(1)\d^{CD} e^{-\f{C\{|\l_3^*|(1+\tau)\}^2}{1+\tau}}
=O(1)\d^{CD} e^{-C(1+\tau)}.
\end{eqnarray*}
\end{itemize}
Keep in mind that each individual wave strength is controlled by the total wave strength
by \eqref{CD} and \eqref{Ri}, that is,
$$
\d^{R_1}+\d^{R_3}+\d^{CD}\leq C\d.
$$
Hence, in summary, it follows from \eqref{Q1}, \eqref{Q2} and the above arguments that
\begin{equation*}
|(Q_{11},Q_{21})|=O(1)\d e^{-C\{|y|+(1+\tau)\}},
\end{equation*}
for some positive constant $C$ independent of $\tau$ and $y$.

\setcounter{equation}{0}
\section{Proof of the main result}

In this section, we shall prove the main result Theorem \ref{limit-th}. By virtue of the arguments in
Section 2.3, it is sufficient to show \eqref{con} besides the regularity of the solution.
To this end, we first reformulate the problem.
\subsection{Reformulation of the problem}
Set the perturbation around the wave profile $(V,U,\T)(\tau,y)$ by
\begin{equation*}
(\phi,\psi,\zeta)(\tau,y)=(v,u,\t)(\tau,y)-(V,U,\T)(\tau,y).
\end{equation*}
Then, after a straightforward calculation, the perturbation $(\phi,\psi,\zeta)(\tau,y)$ satisfies
the system
\begin{equation}
\left\{
\begin{array}{ll}
\di \phi_{\tau}-\psi_{y}=0, \\
\di \psi_{\tau}+(p-P)_{y}=(\frac {u_{y}}v-\f{U_y}{V})_{y}- Q_1,\\[2mm]
\di \f{R}{\g-1}\z_{\tau}+
(pu_y-PU_y)=\nu(\frac{\theta_{y}}{v}-\f{\T_y}{V})_y+(\f{u_y^2}{v}-\frac{U^2_{y}}{V})-Q_2,\\[4mm]
\di (\phi,\psi,\z)(\tau=0,y)=(\phi_0,\psi_0,\z_0)(y),
\end{array}
\right.
 \label{P}
\end{equation}
where the initial data $(\phi_0,\psi_0,\z_0)(y)$ and its derivatives are sufficiently
smooth away from but up to $y=0$, and
$$
(\phi_0,\psi_0,\z_0)(y)\in L^2(\mathbf{R}), \phi_{0y}\in
L^2(\bf R^\pm).
$$
For simplicity, denote
\begin{equation*}
\mathcal{N}_0:=\|(\p_0,\psi_0,\z_0)\|^2+\|\phi_{0y}\dpr^2.
\end{equation*}
In order to prove \eqref{con}, we easily see that it suffices to show
\begin{proposition}\label{prop1} There exists a positive
constant $\d_0$, such that if the wave strength $\d$ and the initial data satisfy
$$
\d+\mathcal{N}_0\leq \d_0,
$$
then the problem \eqref{P} admits a unique global solution
$(\p,\psi,\z)(t,y)$ satisfying
\begin{itemize}
\item[(i)] There exists a positive constant $C$ independent of $t$,
such that
\begin{equation*}
\sup_{\tau\geq0}\Big(\|(\p,\psi,\z)(\tau,\cdot)\|^2+\|\phi_y(\tau,\cdot)\dpr^2\Big)
+\int_0^{+\i}\|(\p_y,\psi_y,\z_y)(\tau,\cdot)\dpr^2d\tau \leq C(\mathcal{N}_0+\d^{\f14}).
\end{equation*}
\item[(ii)] For any $\tau_0>0$, there exists a positive constant $C=C(\tau_0)$, such that
\begin{equation*}
\sup_{\tau\geq\tau_0}\|(\psi_y,\z_y,\psi_\tau,\z_\tau)(\tau,\cdot)\dpr^2
+\int_{\tau_0}^{+\i}\|(\psi_{yy},\z_{yy},\psi_{y\tau},\z_{y\tau})(\tau,\cdot)\dpr^2d\tau
\leq C(\tau_0)(\mathcal{N}_0+\d^{\f14}).
\end{equation*}
\item[(iii)] The jump condition of $\phi(\tau,y)$ at $y=0$ admits the bound
\begin{equation}
|[\phi](\tau)|\leq Ce^{-c\tau}\label{JE}
\end{equation}
where the positive constants $C$ and $c$ are independent of
$\tau\in(0,+\i)$.
\end{itemize}
\end{proposition}

Assume that Proposition \ref{prop1} holds, then for any $\tau_0>0$, one has
\begin{equation*}
\int_{\tau_0}^{+\i}\Big(\|(\phi_y,\psi_y,\z_y)\dpr^2+|\f{d}{d\tau}\|(\phi_{y},\psi_{y},\z_{y})\dpr^2|\Big)d\tau <+\i ,
\end{equation*}
whence,
$$
\lim_{\tau\rightarrow\i}\|(\phi_y,\psi_y,\z_y)\dpr^2=0,
$$
which, together with Proposition \ref{prop1} and Sobolev's inequality, implies that
$$
\lim_{\tau\rightarrow\i}\sup_{y\neq0}\|(\phi,\psi,\z)\|_{L^\i}^2\leq
C\lim_{\tau\rightarrow\i}\|(\phi,\psi,\z)\|\|(\phi_y,\psi_y,\z_y)\dpr\ \leq
C\lim_{\tau\rightarrow\i}\|(\phi_y,\psi_y,\z_y)\dpr\ =0.
$$
The above inequality combined with \eqref{JE} gives \eqref{con}. Thus, the main result Theorem \ref{limit-th} follows from  \eqref{con} and
\eqref{profile-s}.

Denote
\begin{equation*}
\begin{array}{l}
\di N(\tau_*,\tau^*)=\sup_{\tau\in[\tau_*,\tau^*]}\Big\{\|(\p,\psi,\z)(\tau,\cdot)\|^2
+\|(\phi_y,\psi_y,\z_y)(\tau,\cdot)\dpr^2+\|(\psi_\tau,\z_\tau)(\tau,\cdot)\|^2\Big\},\\[4mm]
\di N(\tau_*)=N(\tau_*,\tau_*),
\end{array}
\end{equation*}
and define the solution space by
\begin{equation*}
X[\tau_*,\tau^*]=\left\{(\p,\psi,\z)\left|
\begin{array}{l}
\di (\p,\psi,\z)(\tau,y)\in
 C([\tau_*,\tau^*];H^1({\bf R}^\pm)),\\[1mm]
\di (\psi_y,\z_y)\in L^2(\tau_*,\tau^*; H^1({\bf R}^\pm)),~ \p_y\in L^2(\tau_*,\tau^*;L^2({\bf
 R}^\pm)),\\\di (\psi_\tau,\z_\tau)\in L^\i(\tau_*,\tau^*; L^2({\bf R}^\pm))\cap L^2(\tau_*,\tau^*;H^1({\bf R}^\pm)).
\end{array}
\right.\right\}
\end{equation*}

Since the local existence of solutions to \eqref{P} is
proved in \cite{H3}, we just state it and omit its proof for brevity.

\begin{proposition}\label{local-e}(Local existence) Suppose that
$\mathcal{N}_0$ and the wave strength $\d$ are suitably small such that
$\inf v_0$ and $\inf \t_0$ are positive. Then there exists a
positive time $\tau_0=\tau_0(N(0),\d)>0$, such that the Cauchy problem
\eqref{P} admits a unique solution $(\p,\psi,\z)(\tau,y)\in
X[0,\tau_0]$ satisfying
\begin{equation*}
A(\tau_0)+B(\tau_0)+F(\tau_0)\leq C(\mathcal{N}_0+\d),
\end{equation*}
where
\begin{eqnarray*} &&
A(\tau_0)=\di \sup_{0\leq \tau\leq
\tau_0}\Big\{\|(\p,\psi,\z)(\tau,\cdot)\|^2+\|\phi_y\dpr^2\Big\}+\int_0^{\tau_0}\|(\psi_y,\z_y)\|^2
d\tau,  \\
&& B(\tau_0)=\di \sup_{0\leq \tau\leq
\tau_0}\Big\{g(\tau)^\f12\|\psi_y\|^2+g(\tau)\|\phi_y\dpr^2\Big\}
+\int_0^{\tau_0}g(\tau)^{\f12+\vartheta}(\|\psi_\tau\|^2+\|(\f{u_y}{v})_y\dpr^2)d\tau
\\
&&\di\qquad +\int_0^{\tau_0}g(\tau)(\|\psi_y^2\|^2+\|\t_\tau\|^2+\|(\f{\t_y}{v})_y\dpr^2)d\tau,\\[4mm]
&& F(\tau_0)=\di \sup_{0\leq \tau\leq
\tau_0}\Big\{g(\tau)^{\f32+\vartheta}(\|\psi_\tau\|^2+\|(\f{u_y}{v})_y\dpr^2)
+g(\tau)^3(\|\z_\tau\|^2+\|(\f{\t_y}{v})_y\dpr^2)\Big\}
 \\
&& \qquad \di +\int_0^{\tau_0}g(\tau)^{\f32+\vartheta}\|\psi_{y\tau}\dpr^2+g(\tau)^3\|\z_{y\tau}\dpr^2)d\tau ,
\end{eqnarray*}
with $g(\tau)=\tau\wedge1=\min\{\tau,1\}$ and $\vartheta\in (0,1)$. Moreover,
$v,u,\t$ have the same regularity as in Theorem \ref{limit-th}.
Thus, $v,u_x,\t_x$ have one-side limit at $y=0$ and satisfy the jump conditions
$$
\Big[p-\f{u_y}{v}\Big]=\Big[\f{\t_y}{v}\Big]=0.
$$
Finally, one has the following estimate on the jump at $y=0$,
\begin{equation*}
|[v](\tau)|\leq C\d e^{-c\tau},\qquad \tau >0
\end{equation*}
for some positive constants $C$ and $c$ independent of $\tau$.
\end{proposition}

Hence, in view of the local existence and the
standard continuation process, we see that to prove Proposition \ref{prop1},
it suffices to show the following (uniform) a priori estimate.
\begin{proposition}\label{priori} (A priori estimate) Suppose
that the Cauchy problem \eqref{P} has a solution
$(\p,\psi,\z)(\tau,y)\in X[\tau_1,\tau_2]$. There exists a positive constant
$\eta_1$, such that if
\begin{equation}\label{priori-a}
N(\tau_1,\tau_2)+\d\leq \eta_1,
\end{equation}
 then,
\begin{equation}
N(\tau_1,\tau_2)+\int_{\tau_1}^{\tau_2}\Big\{\|\p_y(\tau,\cdot)\dpr^2
+\|(\psi_y,\z_y)(\tau,\cdot)\dpr_{1}^2+\|(\psi_{y\tau},\z_{y\tau})\dpr^2(\tau,\cdot)\Big\}d\tau
\leq C(N(\tau_1)+\d^{\f14}), \label{p1}
\end{equation}
where the positive constant $C$ is independent of $\tau$.
\end{proposition}

\subsection{Energy estimates}

In this section we will derive the a priori estimate given in Proposition \ref{priori}.
Note that under the a priori assumption \eqref{priori-a}, if $\eta\ll1,$ then if holds that
 $$
\inf_{[\tau_1,\tau_2]\times\mathbf{R}}\{(V+\phi, \T+\z)(\tau,y)\}\geq C_0
 $$
 for some positive constant $C_0$. First, one has the following Lemma:
\begin{lemma}\label{lemma1} Under the assumptions of Proposition \ref{priori},
there exists a constant $C>0$, such that for any $\tau\in [\tau_1,\tau_2]$,
\begin{equation*}
\begin{array}{ll}
\di
\|(\p,\psi,\z,\phi_y)(\tau,\cdot)\dpr^2+\int_{\tau_1}^\tau \Big\{\|\sqrt{(U^{R_1}_y,U^{R_3}_y)}(\p,\z)\|^2+\|(\p_y,\psi_y,\z_y)\dpr^2\Big\} d\tau\\[3mm]
\leq \di  C\|(\p,\psi,\z,\p_{y})\dpr^2(\tau_1) \di
+C\int_{\tau_1}^\tau (1+\tau)^{-\f76}\|(\p,\psi,\z)(\cdot,\tau)\|^2d\tau+C\d^{\f14}\\[4mm]
\di ~~+C\d \int_{\tau_1}^\tau  \dashint_{\bf R}
(1+\tau)^{-1}e^{-\f{c_0y^2}{1+\tau}}|(\p,\z)|^2dy d\tau.
\end{array}
\end{equation*}
\end{lemma}
\noindent{\bf Proof:} Let
$$
\Phi(z)=z-1-\ln z.
$$
Arguing similarly to that in \cite{HLM} or \cite{HWY}, one can get the following equality
\begin{equation}
\begin{array}{ll}
& \di
I_{1\tau}(\tau,y)+H _{1y}(\tau,y)+\f{\T\psi_ y^2}{v\t}+\nu\f{\T\z _y^2}{v\t^2}
+P(U^{R_1}_y+U^{R_3}_y)\left(\Phi(\f{\t V}{v\T})+\g\Phi(\f{v}{V})\right)\\[3mm]
&\di =Q_3-Q_1\psi-Q_2\f{\z}{\t},
\end{array}
\label{(3.13)}
\end{equation}
where
\begin{equation*}
I_1(\tau,y)=R\T\Phi(\f{v}{V})+\f{\psi^2}{2}+\f{R\T}{\g-1}\Phi(\f{\t}{\T}),
\label{I1}
\end{equation*}
\begin{equation}
H_1(\tau,y)=(p-P)\psi-(\f{u _y}{v}-\f{U_ y}{V})\psi-\nu(\f{\t
_y}{v}-\f{\T _y}{V})\f{\z}{\t}, \label{H1}
\end{equation}
and
\begin{equation}
\begin{array}{ll}
Q_3=&\di -PU^{CD}_{y}\left(\Phi(\f{\t
V}{v\T})+\g\Phi(\f{v}{V})\right)+\left(\nu(\f{\T_y}{V})_y+\f{U_y^2}{V}
+ Q_2\right)\Big\{(\g-1)\Phi(\f{v}{V})\\[4mm]
&\di+\Phi(\f{\t}{\T})-\f{\z^2}{\t\T}\Big\} - (\f{1}{v}-\f{1}{V})U_
y\psi_y +(\f1v-\f1V)U_y^2\f{\z}{\t}+2\f{\z\psi_y U_y}{v\t}+\nu
\f{\T_y\z_y\z}{v\t^2}\\[4mm]
&\di -\nu(\f{1}{v}-\f{1}{V})\f{\T\T_y\z_y}{\t^2}+\nu(\f{1}{v}-\f{1}{V})\f{\z\T_y^2}{\t^2} .
\end{array}
\label{Q3}
\end{equation}
Integration of the equality \eqref{(3.13)} with respect to $y$ and $\tau$ over
${\mathbf R}^\pm\times [\tau_1,\tau]$ yields that
\begin{equation}\label{20}
\begin{array}{ll}
\di \int I_1(\tau,y)dy+\int_{\tau_1}^\tau\big[H_1\big](\tau)d\tau+\int_{\tau_1}^\tau\dashint_{\bf R}
\bigg(\f{\T\psi_ y^2}{v\t}+\nu\f{\T\z _y^2}{v\t^2}\bigg)dyd\tau\\[4mm]
\qquad \di +\int_{\tau_1}^\tau\dashint_{\bf R}P(U^{R_1}_y+U^{R_3}_y)\left(\Phi(\f{\t V}{v\T})
+\g\Phi(\f{v}{V})\right)dyd\tau \\ [4mm]
\di=\int I_1(\tau_1,y)dy+\int_{\tau_1}^\tau\dashint_{\bf R} \big(Q_3-Q_1\psi-Q_2\f{\z}{\t}\big)dyd\tau.
\end{array}
\end{equation}

It is easy to observe that the jump of $H_1$ in \eqref{H1} across $y=0$ vanishes, i.e.,
\begin{equation*}
\begin{array}{ll}
\di \big[H_1\big](\tau)&\di=\big[(p-\f{u_y}{v})\psi\big]-\big[(P-\f{U_y}{V})\psi\big]
-\nu\big[(\f{\t_y}{v}-\f{\T_y}{V})\f{\z}{\t}\big]\\
&\di =\big[p-\f{u_y}{v}\big]\psi(\tau,0)-\big[P-\f{U_y}{V}\big]\psi(\tau,0)
-\nu\Big(\big[\f{\t_y}{v}\big]-\big[\f{\T_y}{V}\big]\Big)\f{\z(\tau,0)}{\t(\tau,0)}=0.
\end{array}
\end{equation*}
Recalling that
$$
\Phi(1)=\Phi^\prime(1)=0,\qquad \Phi^{\prime\prime}(z)=z^{-2}>0,
$$
there exists a positive constant $C$, such that if $z$ is near 1, then
$$
C^{-1}(z-1)^2\leq \Phi(z)\leq C(z-1)^2.
$$
Thus under the a priori assumptions \eqref{priori-a}, one gets
\begin{equation}
C^{-1}|\p|^2\leq \Phi(\f{v}{V})\leq C|\p|^2,\qquad C^{-1}|\z|^2\leq
\Phi(\f{\t}{\T})\leq C|\z|^2  \label{(3.17)}
\end{equation}
and
\begin{equation}
C^{-1}|(\p,\z)|^2\leq\Phi(\f{\t V}{v\T})+\g\Phi(\f{v}{V})\leq
C|(\p,\z)|^2. \label{(3.18)}
\end{equation}
Now it follows from \eqref{Q3}, \eqref{(3.17)}, \eqref{(3.18)}
 and Cauchy-Schwarz's inequality that
\begin{equation}
\begin{array}{ll}
\di |Q_3|\leq&\di  \f{\T\psi_ y^2}{4v\t}+\f{\nu\T\z
 _y^2}{4v\t^2}+C\Big\{(|\T^{CD}_y|^2,|\T^{CD}_{yy}|)+(|(V^{R_1}_{y},U^{R_1}_{y},\T^{R_1}_{y})|^2,|\T^{R_1}_{yy}|)\\[3mm]
 &\di +(|(V^{R_3}_{y},U^{R_3}_{y},\T^{R_3}_{y})|^2,|\T^{R_3}_{yy}|)+|Q_2|\Big\}(\p^2+\z^2).
\end{array}
\label{(3.19)}
\end{equation}
By the properties of the viscous contact wave, one can obtain
\begin{equation*}
\int_{\tau_1}^\tau\dashint_{\bf R}
(|\T^{CD}_{y}|^2,|\T^{CD}_{yy}|)(\p^2+\z^2)dyd\tau\leq
   C \d  \int_{\tau_1}^\tau\dashint_{\bf R}(1+\tau)^{-1}e^{-\f{c_0 y^2}{1+\tau}}|(\p,\z)|^2dyd\tau ,
\end{equation*}
while by the properties of the
approximate rarefaction wave in Lemma \ref{lemma-R-NS}, we have that for $i=1,3,$
\begin{equation*}
\begin{array}{ll}
& \di \int_{\tau_1}^\tau\dashint_{\bf R}
(|(V^{R_i}_{y},U^{R_i}_{y},\T^{R_i}_{y})|^2,|\T^{R_i}_{yy}|)(\p^2+\z^2)dyd\tau\\[4mm]
& \di\leq \int_{\tau_1}^\tau
(\|(V^{R_i}_{y},U^{R_i}_{y},\T^{R_i}_{y})\|^2 +\|\T^{R_i}_{yy}\|_{L^1} )\|(\p,\z)\|^2_{L^\infty}d\tau\\[4mm]
&\di \leq C\int_{\tau_1}^\tau
(1+\tau)^ {-1}\|(\p,\z)\|\|(\p_y,\z_y)\|d\tau\\[4mm]
&\di \leq\mu\int_{\tau_1}^\tau  \|(\p_y, \z_y)\|^2d\tau+C_\mu \int_{\tau_1}^\tau (1+\tau)^ {-2}\|(\p,\z)\|^2d\tau ,
\end{array}
\end{equation*}
where and in the sequel $\mu$ is a small positive constant to be determined and $C_\mu$
is some positive constant depending on $\mu$.

Now, it remains to estimate the terms $Q_1\psi$, $Q_2\f{\z}{\t}$ on the right-hand side
of \eqref{20} and the term $|Q_2|(\p^2+\z^2)$ on the right-hand side of \eqref{(3.19)}.
For simplicity, we only estimate $Q_2\f{\z}{\t}$. By \eqref{Q2}, we find that
\begin{equation*}
\begin{array}{ll}
\di \int_{\tau_1}^\tau \dashint_{\mathbf{R}}|Q_{2}\f{\z}{\t}|dy d\tau
\leq C\int_{\tau_1}^\tau  \|\z\|_{L^\i_y}\|Q_2\|_{L^1_y} d\tau\\[4mm]
\quad\di \leq C\int_{\tau_1}^\tau \|\z\|^{\f12}\|\z_y\|^{\f12}\Big(\|Q_{21}\|_{L^1_y}
+\|Q_{22}\|_{L^1_y}+\|Q^{CD}\|_{L^1_y}\Big)d\tau\\[4mm]
\quad\di \leq C\int_{\tau_1}^\tau \|\z\|^{\f12}\|\z_y\|^{\f12}\Big(\d e^{-C(1+\tau)}
+(\d^{r_1}+\d^{r_3})^{\f18}(1+\tau)^{-\f78}+\d (1+\tau)^{-\f32}\Big)d\tau\\[4mm]
\quad\di \leq \mu\int_{\tau_1}^\tau \|\z_y\|^2 d\tau+C_\mu~ \d^{\f16}\int_{\tau_1}^\tau
\|\z\|^{\f23}(1+\tau)^{-\f76}  d\tau\\[4mm]
\quad\di \leq\mu \int_{\tau_1}^\tau \|\z_y\|^2
 d\tau+C_\mu \int_{\tau_1}^\tau\|\z\|^2 (1+\tau)^{-\f76} d\tau+C_\mu~\d^{\f14}.
\end{array}
\end{equation*}
Similarly, one can control the term $Q_1\psi$ and $|Q_2|(\p^2+\z^2)$.

Thus, substituting all the above estimates into \eqref{20} and choosing $\mu$ in the front
of the integral $\di \int_{\tau_1}^\tau\|(\psi_y,\z_y)\|^2d\tau$ small enough,
so that the integral can be absorbed by the left-hand side of \eqref{20}, one concludes
\begin{equation}
\begin{array}{ll}
\di
\|(\p,\psi,\z)(\tau,\cdot)\|^2+\int_{\tau_1}^\tau\big\{\|(\psi_y,\z_y)(\tau,\cdot)\|^2
+\|\sqrt{(U^{R_1}_y, U_y^{R_3})}(\p,\z)(\tau,\cdot)\|^2\big\}d\tau\\
\di \leq C\|(\p,\psi,\z)(\tau_1,\cdot)\|^2+C\int_{\tau_1}^\tau(1+\tau)^{-\f76}\|(\p,\psi,\z)\|^2  d\tau+C\d^{\f14}\\
\di +C\mu\int_{\tau_1}^\tau\|\phi_{y}(\tau,\cdot)\dpr^2d\tau
+C\d \int_{\tau_1}^\tau\dashint_{\bf R}(1+\tau)^{-1}e^{-\f{c_0 y^2}{1+\tau}}|(\p,\z)|^2dy d\tau.
\end{array}\label{(3.25)}
\end{equation}

Next, we estimate $\|\p_y\|^2$. Denote $\tilde v=\f{v}{V}.$
From the system $\eqref{P}_2$, one has
$$
(\f{\tilde{v}_y}{\tilde v})_\tau-\psi_\tau-(p-P)_y-Q_1=0.
$$
Multiplying the above equation by $\f{\tilde{v}_y}{\tilde v}$ and
noticing that
$$
-(p-P)_y=\f{R\t}{v}\f{\tilde{v}_y}{\tilde v}-\f{R\z_y}{v}+(p-P)\f{V_y}{V}-R\T_y(\f1V-\f1v),
$$
one obtains
\begin{eqnarray*} &&
\di \left(\f{1}{2}(\f{\tilde{v}_y}{\tilde v})^2-\psi\f{\tilde{v}_y}{\tilde v}\right)_\tau
+\left(\psi\f{\tilde{v}_\tau}{\tilde v}\right)_y+\f{R\t}{v}(\f{\tilde{v}_y}{\tilde v})^2\\
&& \quad = \di\psi_y(\f{u_y}{v}-\f{U_y}{V})+\left(\f{R\z_y}{v}-(p-P)\f{V_y}{V}
+R\T_y(\f1V-\f1v)-Q_1\right)\f{\tilde{v}_y}{\tilde v}.
\end{eqnarray*}
Integrating the above equality with respect to $y$ and $\tau$ over ${\bf R}^\pm\times[\tau_1,\tau]$
and using Cauchy-Schwarz's inequality, we infer that
\begin{equation}
\begin{array}{ll}
&\di \dashint_{\bf R}\left(\f{1}{2}(\f{\tilde{v}_y}{\tilde
v})^2-\psi\f{\tilde{v}_y}{\tilde v}\right)(\tau,y)dy+\int_{\tau_1}^\tau\left[\psi\f{\tilde{v}_\tau}{\tilde
v}\right](\tau)d\tau+\int_{\tau_1}^\tau\dashint_{\bf
R}\f{R\t}{2v}(\f{\tilde{v}_y}{\tilde v})^2dy d\tau\\[4mm]
\leq&\di \dashint_{\bf R}\left(\f{1}{2}(\f{\tilde{v}_y}{\tilde
v})^2-\psi\f{\tilde{v}_y}{\tilde
v}\right)(\tau_1,y)dy+\int_{\tau_1}^\tau\dashint_{\bf R} |\psi_y(\f{u_y}{v}-\f{U_y}{V})| dy d\tau\\
& \di+C\int_{\tau_1}^\tau\dashint_{\bf R} \left|\f{R\z_y}{v}-(p-P)\f{V_y}{V}
+R\T_y(\f1V-\f1v)-Q_1\right|^2 dy d\tau ,
\end{array}\label{(3.27)}
\end{equation}
where the jump across $y=0$ can be bounded as follows.
\begin{equation*}
\begin{array}{ll}
\di \int_{\tau_1}^\tau\left[\psi\f{\tilde{v}_\tau}{\tilde
v}\right](\tau)d\tau =\int_{\tau_1}^\tau\psi(\tau,0)\left[\f{u_y}{v}-\f{U_y}{V}\right]
(\tau)d\tau=\int_{\tau_1}^\tau\psi(\tau,0)\left[p\right](\tau)d\tau\\[4mm]
\qquad\di =R\int_{\tau_1}^\tau\psi(\tau,0)\t(\tau,0)\left[\f1v\right](\tau)d\tau=-R\int_{\tau_1}^\tau
\f{\psi(\tau,0)\t(\tau,0)}{v(\tau,0+)v(\tau,0-)}\left[v\right](\tau)d\tau\\[4mm]
\qquad\di \leq C\int_{\tau_1}^\tau\|\psi\|_{L^\i}(\tau)|[v]|(\tau_1)e^{-C(\tau-\tau_1)}d\tau
\leq C\d\int_{\tau_1}^\tau\|\psi\|^{\f12}\|\psi_y\|^{\f12}e^{-C(\tau-\tau_1)}d\tau\\[4mm]
\qquad\di \leq \d\int_{\tau_1}^\tau\|\psi_y\|^2d\tau+\d\sup_{\tau\in[\tau_1,\tau_2]}\|\psi\|^2(\tau)+C\d.
\end{array}
\end{equation*}

Using the equality
$$
\f{\tilde{v}_y}{\tilde v}=\f{v_y}{v}-\f{V_y}{V}=\f{\phi_y}{v}-\f{V_y\phi}{vV},
$$
we see that
\begin{equation*}
C^{-1}(|\p_y|^2-|V_y\p|^2)\leq (\f{\tilde{v}_y}{\tilde v})^2\leq
C(|\p_y|^2+|V_y\p|^2).
\end{equation*}
From the definition of $Q_1$ in \eqref{Q1} it follows that
\begin{equation*}
\begin{array}{ll}
\di\int_{\tau_1}^\tau\|Q_1\|^2d\tau\leq C\int_{\tau_1}^\tau\Big(\|Q_{11}\|^2+\|Q_{12}\|^2\Big)d\tau\\
\di \qquad\leq C\int_{\tau_1}^\tau\Big(\|Q_{11}\|^2+\|(U^{R_1}_{yy},U^{R_3}_{yy},U^{R_1}_yV^{R_1}_y,U^{R_3}_yV^{R_3}_y)\|^2\Big)d\tau\leq C\d^{\f14}.
\end{array}
\end{equation*}

Therefore, substituting all the above estimates into \eqref{(3.27)}, we conclude that
\begin{equation}
\begin{array}{ll}
\di \quad \|\phi_y(\tau,\cdot)\dpr^2+\int_{\tau_1}^\tau\|\p_y\dpr^2 d\tau \leq
C\|(\p,\psi,\p_{y})\dpr^2(\tau_1)+C\|(\p,\psi)(\tau,\cdot)\|^2\\
\di +C\d \int_{\tau_1}^\tau\dashint_{\bf R}(1+\tau)^{-1}e^{-\f{c_0 y^2}{1+\tau}}|(\p,\z)|^2dy
d\tau+C\int_{\tau_1}^\tau\|(\psi_y,\z_y)\|^2d\tau
\\[4mm]
\di+C\int_{\tau_1}^\tau(1+\tau)^{-\f76}\|(\p,\psi,\z)\|^2 d\tau+C\d^{\f14}.
\end{array}
 \label{(3.30)}
\end{equation}

Multiplying the inequality \eqref{(3.25)} by a large constant
$C_1>0$, and summing the resulting inequality with \eqref{(3.30)}, we obtain
Lemma \ref{lemma1}. This completes the proof. $\hfill\Box$
\vspace{2mm}

Next, we derive the higher order estimates, which are summarized in the following Lemma:

\begin{lemma}\label{lemma2} Under the assumptions of Proposition \ref{priori}, it holds that
\begin{equation*}
\begin{array}{ll}
\di N(\tau_1,\tau_2)+\int_{\tau_1}^{\tau_2}
\big\{\|\sqrt{(U^{R_1}_y,U^{R_3}_y)}(\p,\z)\|^2+\|\p_y\dpr^2+\|(\psi_y,\z_y)\dpr_1^2
+\|(\psi_{y\tau},\z_{y\tau})\dpr^2\big\} d\tau\\
\di \leq CN(\tau_1)+C\int_{\tau_1}^{\tau_2}(1+\tau)^{-\f76}\|(\p,\psi,\z)\|^2
d\tau+C\d^{\f14}+C \d \int_{\tau_1}^{\tau_2} \dashint_{\bf R}
(1+\tau)^{-1}e^{-\f{c_0 y^2}{1+\tau}} |(\p,\z)|^2dy d\tau.
\end{array}
\end{equation*}
\end{lemma}
\noindent{\rm\bf Proof:} Multiplying the equation
$\eqref{P}_2$ by $\di-\psi_{yy}$, one gets
\begin{equation*}
\begin{array}{ll}
\di \left(\f{\psi_y^2}{2}\right)_\tau-\left(\psi_\tau\psi_y\right)_y+\f{\psi_{yy}^2}{v}
=\Big\{(p-P)_y+\f{v_y}{v^2}\psi_y-\big(U_y (\f{1}{v}-\f{1}{V})\big)_y+Q_1\Big\}\psi_{yy}.
\end{array}
\end{equation*}
Integration of the above equation with respect to $y$ and $\tau$ over ${\bf R}^\pm\times[\tau_1,\tau]$ gives
\begin{equation}
\begin{array}{ll}
\di \dashint_{\bf R}\f{\psi_y^2}{2}(\tau,y)dy+\int_{\tau_1}^\tau\dashint_{\bf R}\f{\psi_{yy}^2}{v}dyd\tau
=\dashint_{\bf R}\f{\psi_y^2}{2}(\tau_1,y)dy-\int_{\tau_1}^\tau\left[\psi_\tau\psi_y\right](\tau)d\tau\\
\di ~~+\int_{\tau_1}^\tau\dashint_{\bf R}\Big\{(p-P)_y+\f{v_y}{v^2}\psi_y
-\big(U_y (\f{1}{v}-\f{1}{V})\big)_y+Q_1\Big\}\psi_{yy}dyd\tau =:\sum_{i=1}^3J_i.
\end{array}\label{(3.31)}
\end{equation}

We have to estimate $J_i$. First, the jump $J_2$ can be bounded as follows.
\begin{equation}\label{J2-1}
\begin{array}{ll}
J_2&\di =-\int_{\tau_1}^\tau\left[\psi_\tau\psi_y\right](\tau)d\tau=-\int_{\tau_1}^\tau\psi_\tau(\tau,0)\left[\psi_y\right](\tau)d\tau\\[3mm]
 &\di= -\int_{\tau_1}^\tau\psi_\tau(\tau,0)\left[u_y\right](\tau)d\tau=-\int_{\tau_1}^\tau\psi_\tau(\tau,0)\left[(\f{u_y}{v}-p)v\right](\tau)d\tau\\[3mm]
&\di=-\int_{\tau_1}^\tau\psi_\tau(\tau,0)(\f{u_y}{v}-p)(\tau,0)\left[v\right](\tau)d\tau\\[3mm]
&\di\leq C\int_{\tau_1}^\tau\|\psi_\tau\|_{L^\i}\big(\|\psi_y\|_{L^\i}+1\big)\left[v\right](\tau_1)e^{-C(\tau-\tau_1)}d\tau\\[3mm]
&\di\leq C\d\int_{\tau_1}^\tau\|\psi_\tau\|^{\f12}\|\psi_{y\tau}\dpr^{\f12}\big(\|\psi_y\|^{\f12}\|\psi_{yy}\dpr^{\f12}+1\big)e^{-C(\tau-\tau_1)}d\tau.
\end{array}
\end{equation}
In view of $\eqref{P}_2$ and \eqref{priori-a}, one has
\begin{equation}\label{psi-tau}
\begin{array}{ll}
\|\psi_\tau\|&\di \leq C\Big(\|\psi_{yy}\dpr+\|(\phi_y,\psi_y,\z_y)\dpr+\|(U_{yy},V_y,U_y,\T_y)\phi\|+\|Q_1\|\Big)\\[3mm]
&\di \leq C\Big(\|\psi_{yy}\dpr+\|(\phi_y,\psi_y,\z_y)\dpr+\d\Big).
\end{array}
\end{equation}
Substituting \eqref{psi-tau} into \eqref{J2-1}, we obtain
\begin{equation}\label{J2}
\begin{array}{ll}
\di |J_2| \leq C\d\int_{\tau_1}^\tau\Big(\|\psi_{yy}\dpr+\|(\phi_y,\psi_y,\z_y)\dpr+\d\Big)^{\f12}
\|\psi_{y\tau}\dpr^{\f12}\Big(\|\psi_y\|^{\f12}\|\psi_{yy}\dpr^{\f12}+1\Big)e^{-C(\tau-\tau_1)}d\tau\\
\di \leq \mu\int_{\tau_1}^\tau\|(\psi_{yy},\psi_{y\tau})\dpr^2d\tau
+C_\mu~\d\int_{\tau_1}^\tau \|(\p_y,\psi_y,\z_y)\dpr^2 d\tau+C_\mu\d .
\end{array}
\end{equation}
On the other hand, $J_3$ can be estimates as follows.
\begin{equation}\label{J3}
\begin{array}{ll}
J_3&\di =\int_{\tau_1}^\tau\dashint_{\bf R}\bigg\{(p-P)_y+\f{v_y}{v^2}\psi_y
-\big(U_y (\f{1}{v}-\f{1}{V})\big)_y+Q_1\bigg\}\psi_{yy}dyd\tau\\[4mm]
&\di \leq C\int_{\tau_1}^\tau\dashint_{\bf R}\Big\{|(\phi_y,\z_y)|+|(\phi,\z)||(\phi_y,V_y,\T_y,U_{yy})|\\
&\di\qquad\qquad\qquad~~+|(\phi_y,V_y)||(\psi_y,U_y,U_y\phi)|+|Q_1|\Big\}|\psi_{yy}|dyd\tau\\[3mm]
&\di \leq \mu \int_{\tau_1}^\tau \|\psi_{yy}\dpr^2d\tau+C_\mu \int_{\tau_1}^\tau
\|(\phi_y,\psi_y,\z_y)\dpr^2d\tau+C_\mu~\int_{\tau_1}^\tau(1+\tau)^{-\f76}\|(\phi,\psi,\z)\|^2d\tau\\[3mm]
&\di ~~~+C_\mu~\d+C_\mu ~\d \int_{\tau_1}^\tau\dashint_{\bf
R}(1+\tau)^{-1}e^{-\f{c_0 y^2}{1+\tau}}|(\p,\z)|^2dy
d\tau.
\end{array}
\end{equation}
Substituting \eqref{J2} and \eqref{J3} into \eqref{(3.31)} and choosing $\mu$ suitably small
in the front of the integral $\int_{\tau_1}^\tau \|\psi_{yy}\dpr^2d\tau$, we deduce that
\begin{equation}
\begin{array}{ll}
\di \|\psi_y\|^2(\tau) +\int_{\tau_1}^\tau\|\psi_{yy}\dpr^2d\tau \leq
C\|\psi_{y}\|^2(\tau_1)+C\mu \int_{\tau_1}^\tau \|\psi_{y\tau}\dpr^2d\tau\\[4mm]
\di\quad +C_\mu\int_{\tau_1}^\tau(1+\tau)^{-\f76}\|(\p,\z)\|^2
d\tau+C_\mu~\d+C_\mu\int_{\tau_1}^\tau \|(\p_y,\psi_y,\z_y)\dpr^2d\tau \\[4mm]
\di \quad +C_\mu~
\d \int_{\tau_1}^\tau \dashint_{\bf R}
(1+\tau)^{-1}e^{-\f{c_0 y^2}{1+\tau}} |(\p,\z)|^2dy d\tau .
\end{array}\label{(3.32)}
\end{equation}

Multiplication of the equation $\eqref{P}_3$ with $-\z_{yy}$ yields that
\begin{equation*}
\begin{array}{ll}
\di \f{R}{\g-1}\left(\f{\z_y^2}{2}\right)_\tau-\f{R}{\g-1}\left(\z_\tau\z_y\right)_y+\nu\f{\z_{yy}^2}{v}\\
\di =\bigg\{(pu_y-PU_y)
+\nu\f{\z_yv_y}{v^2}-\nu\big(\T_y(\f1v-\f1V)\big)_y-(\f{u_y^2}{v}-\f{U_y^2}{V})+Q_2\bigg\}\z_{yy}.
\end{array}
\end{equation*}
Integrating the above equality with respect to $y$ and $\tau$ over ${\bf R}^\pm\times[\tau_1,\tau]$,
and employing almost the same arguments as those used for $\|\psi_y\dpr^2(\tau)$ in \eqref{(3.32)},
we obtain
\begin{equation}
\begin{array}{ll}
&\di \|\z_y\|^2(\tau) +\int_{\tau_1}^\tau\|\z_{yy}\dpr^2d\tau \leq
C\|\z_{y}\|^2(\tau_1)+C\d^{\f14}+C\int_{\tau_1}^\tau(1+\tau)^{-\f76}\|(\p,\z)\|^2
d\tau \\
&\di \quad +C\int_{\tau_1}^\tau \|(\p_y,\psi_y,\z_y)\dpr^2d\tau+C
(\d )^2\int_{\tau_1}^\tau \dashint_{\bf R}
(1+\tau)^{-1}e^{-\f{c_0 y^2}{1+\tau}} |(\p,\z)|^2dy d\tau,
\end{array}\label{(3.33)}
\end{equation}
where we have used the following jump estimate across $y=0$
\begin{equation*}
\begin{array}{ll}
\di -\f{R}{\g-1}\int_{\tau_1}^\tau\left[\z_\tau\z_y\right](\tau) d\tau
=-\f{R}{\g-1}\int_{\tau_1}^\tau\z_\tau(\tau,0)\left[\z_y\right](\tau) d\tau\\[4mm]
\di =-\f{R}{\g-1}\int_{\tau_1}^\tau\z_\tau(\tau,0)\left[\t_y\right](\tau) d\tau
=-\f{R}{\g-1}\int_{\tau_1}^\tau\z_\tau(\tau,0)\f{\t_y}{v}(\tau,0)\left[v\right](\tau) d\tau\\
\di \leq C\int_{\tau_1}^\tau\|\z_\tau\|_{L^\i}\big(1+\|\z_y\|_{L^\i}\big)[v](\tau_1)e^{-C(\tau-\tau_1)}d\tau\\
\di \leq C\d\int_{\tau_1}^\tau\|\z_\tau\dpr^{\f12}\|\z_{y\tau}\dpr^{\f12}
\big(1+\|\z_y\|^{\f12}\|\z_{yy}\dpr^{\f12}\big)e^{-C(\tau-\tau_1)}d\tau
\end{array}
\end{equation*}
and the estimate
\begin{equation}\label{zeta-tau}
\begin{array}{ll}
\di \|\z_\tau\|\leq C\Big(\|\z_{yy}\dpr+\|(\phi_y,\psi_y,\z_y)\dpr+\|(U_y,\T_{yy},\T_yV_y,U_y^2)(\phi,\z)\|+\|Q_2\|\Big)\\
\di \qquad \leq C\Big(\|\z_{yy}\dpr+\|(\phi_y,\psi_y,\z_y)\dpr+\d\Big).
\end{array}
\end{equation}

It follows from \eqref{psi-tau} and \eqref{zeta-tau} that
\begin{equation}\label{tau}
\begin{array}{ll}
\di \int_{\tau_1}^{\tau_2}\|(\psi_\tau,\z_\tau)(\tau,\cdot)\|^2d\tau\\
\di \leq C\Big(\int_{\tau_1}^{\tau_2}\|(\psi_{yy},\z_{yy})\dpr^2d\tau+\int_{\tau_1}^{\tau_2}\|(\phi_y,\psi_y,\z_y)\dpr^2d\tau\\[4mm]
\di\qquad\qquad+\int_{\tau_1}^{\tau_2}\|(U_y,\T_{yy},\T_yV_y,U_y^2)(\phi,\z)\|^2d\tau+\int_{\tau_1}^{\tau_2}\|Q_2\|^2d\tau\Big)\\[4mm]
\di \leq C\int_{\tau_1}^{\tau_2}\|(\psi_{yy},\z_{yy})\dpr^2d\tau+C\int_{\tau_1}^{\tau_2}
\|(\phi_y,\psi_y,\z_y)\dpr^2d\tau+\int_{\tau_1}^{\tau_2}(1+\tau)^{-2}\|(\phi,\z)\|^2d\tau+C\d^{\f14}.
\end{array}
\end{equation}

Now we turn to control $\di \sup_{\tau\in[\tau_1,\tau_2]} \|(\psi_{\tau},\z_{\tau})\dpr^2$.
First, applying the operator $\partial_\tau$ to the equation $(\ref{P})_2$, we get
\begin{equation*}
\psi_{\tau\tau}=\big(\f{u_y}{v}-p\big)_{y\tau}-\big(\f{U_y}{V}-P\big)_{y\tau}-Q_{1\tau}.
\end{equation*}
Multiplication of the above equation by $\psi_\tau$ gives
\begin{equation*}
\begin{array}{ll}
\di \left(\f{\psi_\tau^2}{2}\right)_\tau+\f{\psi_{y\tau}^2}{v}=\Big\{\psi_\tau\big(\f{u_y}{v}
-p\big)_\tau-\psi_\tau\big(\f{U_y}{V}-P\big)_\tau\Big\}_y\\[4mm]
\di\qquad -\psi_{y\tau}\f{U_{y\tau}}{v}+\psi_{y\tau}\f{u_y}{v^2}v_\tau+\psi_{y\tau}(\f{U_y}{V})_\tau
+\psi_{y\tau}(p-P)_\tau-\psi_\tau Q_{1\tau}.
\end{array}
\end{equation*}
If we integrate the above equality with respect to $y$ and $\tau$ over ${\bf R}^\pm\times[\tau_1,\tau]$,
we find that
\begin{equation}\label{f1}
\begin{array}{ll}
\di \dashint_{\bf R}\f{\psi_\tau^2}{2}(\tau,y)dy+\int_{\tau_1}^\tau\dashint_{\bf R}\f{\psi_{y\tau}^2}{v}dyd\tau\\[4mm]
\di~~=\dashint_{\bf R}\f{\psi_\tau^2}{2}(\tau_1,y)dy-\int_{\tau_1}^\tau\Big[\psi_\tau
\big(\f{u_y}{v}-p\big)_\tau-\psi_\tau\big(\f{U_y}{V}-P\big)_\tau\Big](\tau)d\tau\\[4mm]
\di +\int_{\tau_1}^\tau\dashint_{\bf R}\Big\{-\psi_{y\tau}\f{U_{y\tau}}{v}+\psi_{y\tau}
\f{u_y}{v^2}v_\tau+\psi_{y\tau}(\f{U_y}{V})_\tau+\psi_{y\tau}(p-P)_\tau-\psi_\tau Q_{1\tau}\Big\}dyd\tau ,
\end{array}
\end{equation}
where the jump across $y=0$ in fact vanishes, i.e.,
\begin{equation}
\begin{array}{ll}
\di \Big[\psi_\tau\big(\f{u_y}{v}-p\big)_\tau-\psi_\tau\big(\f{U_y}{V}-P\big)_\tau\Big](\tau)\\[3mm]
\di=[\psi_\tau](\tau)\big(\f{u_y}{v}-p\big)_\tau(\tau,0-) +\psi_\tau(\tau,0+)
\Big[\big(\f{u_y}{v}-p\big)_\tau\Big](\tau)-[\psi_\tau](\tau)\big(\f{U_y}{V}-P\big)_\tau(\tau,0)\\[3mm]
\di=[\psi]_\tau(\tau)\big(\f{u_y}{v}-p\big)_\tau(\tau,0-)+\psi_\tau(\tau,0+)
\Big[\f{u_y}{v}-p\Big]_\tau(\tau)-[\psi]_\tau(\tau)\big(\f{U_y}{V}-P\big)_\tau(\tau,0)\\[2mm]
\di =0.  \end{array}   \label{js1}
\end{equation}

Now we apply $\partial_\tau$ to the equation $(\ref{P})_3$ to deduce that
\begin{equation*}
\f{R}{\g-1}\z_{\tau\tau}=\nu\big(\f{\t_y}{v}\big)_{y\tau}-\nu\big(\f{\T_y}{V}\big)_{y\tau}
+\Big\{u_y\big(\f{u_y}{v}-p\big)\Big\}_{\tau}-\Big\{u_y\big(\f{U_y}{V}-P\big)\Big\}_{\tau}-Q_{2\tau}.
\end{equation*}
Multiplying the above equation by $\z_\tau$, one has
\begin{equation*}
\begin{array}{ll}
\di \f{R}{\g-1}(\f{\z_\tau^2}{2})_\tau+\nu\f{\z_{y\tau}^2}{v}
=\Big\{\nu\z_\tau\big(\f{\t_y}{v}\big)_\tau-\nu\z_\tau\big(\f{\T_y}{V}\big)_\tau\Big\}_y\\[3mm]
\di \qquad+\nu\z_{y\tau}\f{\T_{y\tau}}{v}
+\nu\z_{y\tau}\f{\t_y}{v^2}v_\tau+\nu\z_{y\tau}(\f{\T_y}{V})_\tau+\z_\tau u_{y\tau}(\f{u_y}{v}-p)\\[3mm]
\di\qquad+\z_\tau u_y(\f{u_y}{v}-p)_{\tau} -\z_\tau U_{y\tau}(\f{U_y}{V}-P)-\z_\tau U_y(\f{U_y}{V}-P)_{\tau}-\z_\tau Q_{2\tau}.
\end{array}
\end{equation*}
Integrating the above equality with respect to $y$ and $\tau$ over ${\bf R}^\pm\times[\tau_1,\tau]$, we deduce that
\begin{equation}\label{f2}
\begin{array}{ll}
\di \dashint_{\bf R}\f{R\z_\tau^2}{2(\g-1)}(\tau,y)dy+\int_{\tau_1}^\tau\dashint_{\bf R}\nu\f{\z_{y\tau}^2}{v}dyd\tau
=\dashint_{\bf R}\f{R\z_\tau^2}{2(\g-1)}(\tau_1,y)dy\\[5mm]
\di -\int_{\tau_1}^\tau\Big[\nu\z_\tau\big(\f{\t_y}{v}\big)_\tau-\nu\z_\tau\big(\f{\T_y}{V}\big)_\tau\Big](\tau)d\tau
+\int_{\tau_1}^\tau\dashint_{\bf R}\Big\{\nu\z_{y\tau}\f{\T_{y\tau}}{v}+\nu\z_{y\tau}\f{\t_y}{v^2}v_\tau+\nu\z_{y\tau}(\f{\T_y}{V})_\tau\\[5mm]
\di +\z_\tau u_{y\tau}(\f{u_y}{v}-p)+\z_\tau u_y(\f{u_y}{v}-p)_{\tau} -\z_\tau U_{y\tau}(\f{U_y}{V}-P)
-\z_\tau U_y(\f{U_y}{V}-P)_{\tau}-\z_\tau Q_{2\tau}\Big\}dyd\tau ,
\end{array}
\end{equation}
where the jump in fact vanishes.
\begin{equation}
\begin{array}{ll}
\di \Big[\nu\z_\tau\big(\f{\t_y}{v}\big)_\tau-\nu\z_\tau\big(\f{\T_y}{V}\big)_\tau\Big](\tau)\\[2mm]
\di=\nu[\z_\tau](\tau)\big(\f{\t_y}{v}\big)_\tau(\tau,0-)+\nu \z_\tau(\tau,0+)\Big[\big(\f{\t_y}{v}\big)_\tau\Big](\tau)
-\nu[\z_\tau](\tau)\big(\f{\T_y}{V}\big)_\tau(\tau,0)\\[3mm]
\di=\nu[\z]_\tau(\tau)\big(\f{\t_y}{v}\big)_\tau(\tau,0-)+\nu\z_\tau(\tau,0+)\Big[\f{\t_y}{v}\Big]_\tau(\tau)
-\nu[\z]_\tau(\tau)\big(\f{\T_y}{V}\big)_\tau(\tau,0)\\[2mm]
\di=0.
\end{array}   \label{js2}
\end{equation}
Hence, taking into account \eqref{js1} and \eqref{js2}, we get from \eqref{f1} and \eqref{f2} that
\begin{equation}\label{f3}
\begin{array}{ll}
\di \|(\psi_\tau,\z_\tau)\dpr^2(t) +\int_{\tau_1}^\tau\|(\psi_{y\tau},\z_{y\tau})\dpr^2d\tau \leq
C\|(\psi_\tau,\z_\tau)\dpr^2(\tau_1)+C \int_{\tau_1}^\tau \|(\psi_{\tau},\z_\tau)\dpr^2d\tau\\[4mm]
\di\quad +C\int_{\tau_1}^\tau(1+\tau)^{-\f76}\|(\p,\z)\|^2
d\tau+C~\d+C\int_{\tau_1}^\tau \|(\phi_y,\psi_y,\z_y)\|^2d\tau \\[4mm]
\di \quad +C\d \int_{\tau_1}^\tau \dashint_{\bf R}
(1+\tau)^{-1}e^{-\f{c_0 y^2}{1+\tau}} |(\p,\z)|^2dy d\tau.
\end{array}
\end{equation}

Combing the estimates \eqref{(3.32)}, \eqref{(3.33)}, \eqref{tau}, \eqref{f3} and Lemma \ref{lemma1}
together, we obtain Lemma \ref{lemma2}, and  the proof is completed. $\hfill\Box$
\vspace{2mm}

It remains to control the term
$$
\d \int_{\tau_1}^\tau \dashint_{\bf R}
(1+\tau)^{-1}e^{-\f{c_0 y^2}{1+\tau}} |(\p,\z)|^2dy d\tau,
$$
which comes from the viscous contact wave. We shall use the estimate on the heat kernel in
\cite{HLM} to get the desired estimates.

\begin{lemma}\label{lemma3}  Suppose that $Z(t,y)$ satisfies
$$
Z\in L^\i(0,T; L^2(\mathbf{R}^{\pm})),~~Z_y\in L^2(0,T;
L^2(\mathbf{R}^{\pm})),~~Z_\tau\in L^2(0,T; H^{-1}(\mathbf{R}^{\pm})),  $$
then
\begin{equation}
\begin{array}{ll}
&\di \int_{\tau_1}^\tau \dashint_{\mathbf{R}}(1+\tau)^{-1}Z^2e^{-\f{2\b y^2}{1+\tau}}dy d\tau\\[4mm]
&\di \leq C_\b\bigg\{ \|Z(\tau_1,y)\|^2+\int_{\tau_1}^\tau \|h_y\dpr^2
d\tau+\int_{\tau_1}^\tau\langle Z_\tau,Zg_\b^2\rangle_{H^1\times
H^{-1}(\mathbf{R}^{\pm})}d\tau\bigg\}
\end{array}\label{(3.35)}
\end{equation}
where
\begin{equation}\label{g-b}
g_\b(\tau,y)=
\di (1+\tau)^{-\f12}\int_{0}^{y} e^{-\f{\b
\eta^2}{1+\tau}}d\eta
\end{equation}
and $\b>0$ is the constant to be determined. \end{lemma}

\begin{remark}
Lemma \ref{lemma3} can be shown using arguments similar to those in
\cite{HLM}, and hence its proof will be omitted here for simplicity. Note that the domain
considered here consists of two half lines $\mathbf{R}^\pm$, and hence the jump across $y=0$ should
be treated. In view of this, the functional $g_\b$ should be chosen in \eqref{g-b},
so that $g_\b$ is continuous at $y=0$. Furthermore, it holds that $g_\b(\tau,0)\equiv0.$
\end{remark}
\begin{lemma}\label{lemma4}  Under the assumptions of Proposition \ref{priori}, it holds that
\begin{equation*}
\begin{array}{ll}
&\di \int_{\tau_1}^\tau\dashint_{\bf R}
\f{e^{-\f{c_0 y^2}{1+\tau}}}{1+\tau} |(\p,\psi,\z)|^2 dy
d\tau\leq
C\d+C\|(\p,\psi,\z)(\tau_1,\cdot)\|^2+C\|(\p,\psi,\z)(\tau,\cdot)\|^2\\[2mm]
&~~~~~~~~~~~~~~~~~~~~\di +C\int_{\tau_1}^\tau\|(\p_y,\psi_y,\z_y)\dpr^2d\tau+C\int_{\tau_1}^\tau
(1+\tau)^{-\f76}\|(\p,\psi)\|^2 d\tau.
\end{array}\label{(3.36)}
\end{equation*}\end{lemma}
\noindent{\bf Proof:}
From the equation $\eqref{P}_2$ and the fact $p-P=\f{R\z-P\p}{v}$ one gets
\begin{equation}
\psi_\tau+(\f{R\z-P\p}{v})_y=(\f{u_y}{v}-\f{U_y}{V})_y-Q_1.\label{(3.37)}
\end{equation}
Let
\begin{equation*}\label{G-a}
G_\a(\tau,y)=(1+\tau)^{-1}\int_{0}^{y} e^{-\f{\a
\eta^2}{1+\tau}}d\eta ,
\end{equation*}
where $\a$ is a positive constant to be determined.
Multiplying the equation \eqref{(3.37)} by $G_\a(R\z-P\p)$, we find that
\begin{equation}
\begin{array}{ll}
&\di
\left(\f{G_\a(R\z-P\p)^2}{2v}\right)_y-(G_\a)_y\f{(R\z-P\p)^2}{2v}+\f{G_\a(R\z-P\p)^2}{2v^2}(V_y+\p_y)\\[2mm]
&\di =-G_\a (R\z-P\p)\psi_\tau+ G_\a (R\z-P\p)(\f{u_y}{v}-\f{U_y}{V})_y-G_\a
(R\z-P\p)Q_1.
\end{array}\label{(3.38)}
\end{equation}
Noticing that
\begin{equation}
-G_\a (R\z-P\p)\psi_\tau=-\big(G_\a
(R\z-P\p)\psi\big)_\tau +(G_\a)_\tau(R\z-P\p)\psi+G_\a \psi(R\z-P\p)_\tau
\label{(3.39)}
\end{equation}
and
\begin{equation}
\begin{array}{ll}
\di
(R\z-P\p)_\tau=R\z_\tau-P_\tau\p-P\p_\tau\\[2mm]
\di =-\g P\psi_y+(\g-1)\Big\{-(p-P)(U_y+\psi_y)
+(\f{u_y^2}{v}-\f{U_y^2}{V})+\nu(\f{\t_y}{v}-\f{\T_y}{V})_y -Q_2\Big\} -P_\tau\p ,
\end{array}\label{(3.40)}
\end{equation}
if we insert \eqref{(3.40)} into \eqref{(3.39)} and use the equality
\begin{equation*}
\di -G_\a \g P \psi_y \psi =-\Big(\g G_\a P\f{\psi^2}{2}\Big)_y+\g
P(G_\a)_y\f{\psi^2}{2}+\g P_y\f{\psi^2}{2}, \label{(3.41)}
\end{equation*}
we get from \eqref{(3.38)} that
\begin{equation}
\f{e^{-\f{\a y^2}{1+\tau}}}{2(1+\tau)}\Big\{(R\z-P\p)^2+\g
P\psi^2\Big\}=\big\{G_\a v(R\z-P\p)\psi\big\}_\tau+H_{2y}+Q_4, \label{(3.42)}
\end{equation}
where
\begin{equation*}
\di H_2=\f{G_\a(R\z-P\p)^2}{2v}+\g G_\a P\f{\psi^2}{2}-\nu(\g-1)G_\a \psi
(\f{\t_y}{v}-\f{\T_y}{V})-G_\a (R\z-P \p)(\f{u_y}{v}-\f{U_y}{V})
 \label{(3.43)}
\end{equation*}
and
\begin{equation*}
\begin{array}{ll}
\di Q_4= \f{G_\a(R\z-P\p)^2}{2v^2}(V_y+\p_y)-(G_\a)_\tau (R\z-P\p)\psi
+\big(G_\a (R\z-P\p)\big)_y(\f{u_y}{v}-\f{U_y}{V}) \\[2mm]
\di \qquad+(\g-1)G_\a
\psi\left\{(p-P)(U_y+\psi_y)-(\f{u_y^2}{v}-\f{U_y^2}{V})+Q_2\right\}\\[2mm]
\di\qquad +(\g-1)\nu(G_\a \psi)_y(\f{\t_y}{v}-\f{\T_y}{V})+G_\a (R\z-P\p)Q_1+G_\a \psi P_\tau\p-\g P_y\f{\psi^2}{2}.
\end{array}
 \label{(3.44)}
\end{equation*}

Integrating \eqref{(3.42)} over $\mathbf{R}^\pm\times [\tau_1,\tau]$, one infers that
\begin{equation}
\begin{array}{ll}
&\di \int_{\tau_1}^\tau\dashint_{\bf
R}\f{e^{-\f{\a y^2}{1+\tau}}}{1+\tau}\big\{(R\z-P\p)^2+\psi^2\big\}dy d\tau=
\dashint_{\bf R}\big\{G_\a v(R\z-P\p)\psi\big\}(\tau,y)dy\\[4mm]
&\di\qquad -\dashint_{\bf
R}\big\{G_\a v(R\z-P\p)\psi\big\}(\tau_1,y)dy+\int_{\tau_1}^\tau\left[H_{2}\right](\tau)d\tau
+\int_{\tau_1}^\tau\dashint_{\bf R}Q_4dyd\tau.
\end{array}\label{(3.46)-1}
\end{equation}
Here we only analyze the jump term $[H_2]$ across $y=0$, the other terms in \eqref{(3.46)-1}
can be estimated similarly to those in \cite{HLM} or \cite{HWY}. Recalling that $G_\a(\tau,y)$ is continuous
at $y=0$ and $G_\a(\tau,0)\equiv 0$, we easily see that
$$
[H_2](\tau)=0.
$$
Thus, from \eqref{(3.46)-1} one gets
\begin{equation}
\begin{array}{ll}
&\di \int_{\tau_1}^\tau\dashint_{\bf
R}\f{e^{-\f{\a y^2}{1+\tau}}}{1+\tau}\big\{(R\z-P\p)^2+\psi^2\big\}dy d\tau\leq
C\d+C\|(\p,\psi,\z)(\tau_1,\cdot)\|^2\\[3mm]
&\di \qquad\qquad +C\|(\p,\psi,\z)(\tau,\cdot)\|^2+C\int_{\tau_1}^\tau
(1+\tau)^{-\f76}\|(\p,\psi,\z)(\tau,\cdot)\|^2d\tau\\[2mm]
&\di\qquad\qquad +C\int_{\tau_1}^\tau
\|(\p_y,\psi_y,\z_y)(\tau,\cdot)\dpr^2 d\tau+C\d \int_{\tau_1}^\tau\dashint_{\bf
R}\f{e^{-\f{\a y^2}{1+\tau}}}{1+\tau}|(\p,\z)|^2dy
d\tau.
\end{array}\label{(3.46)}
\end{equation}

In order to get the desired estimate in Lemma \ref{lemma4}, we will use Lemma \ref{lemma3} to derive
another similar estimate from the energy equation $\eqref{P}_3$. To this end, we set
$$
Z=\f{R}{\g-1}\z+P\p
$$
in Lemma \ref{lemma3}. Thus we only need to compute the last term in
\eqref{(3.35)}. From the energy equation $\eqref{P}_3$, we have
\begin{equation*}
Z_\tau=P_\tau\p-(p-P)u_y+\nu(\f{\t_y}{v}-\f{\T_y}{V})_y+(\f{u_y^2}{v}-\f{U_y^2}{V})-Q_2,\label{(3.47)}
\end{equation*}
whence
\begin{equation*}
\begin{array}{ll}
&\di \int_{\tau_1}^\tau\langle Z_\tau, Zg_\b^2\rangle_{H^1\times
H^{-1}({\bf R}^\pm)}d\tau\\[5mm]
=&\di\int_{\tau_1}^\tau\dashint_{\mathbf{R}}\big(P_\tau\p-(p-P)U_y\big)Zg_\b^2dy
d\tau-\int_{\tau_1}^\tau\dashint_{\mathbf{R}}(p-P)\psi_y Zg_\b^2dy d\tau\\[4mm]
&\di+\int_{\tau_1}^\tau\Big[\nu(\f{\t_y}{v}-\f{\T_y}{V}) Zg_\b^2\Big](\tau) d\tau-
\int_{\tau_1}^\tau\dashint_{\mathbf{R}}\nu(\f{\t_y}{v}-\f{\T_y}{V})(Zg_\b^2)_y
dy d\tau\\[3mm]
&\di +\int_{\tau_1}^\tau\dashint_{\mathbf{R}}(\f{u^2_y}{v}-\f{U^2_y}{V})Zg_\b^2
dy d\tau-\int_{\tau_1}^\tau\dashint_{\mathbf{R}}Q_2 Zg_\b^2 dy d\tau =:\di \sum_{i=1}^6 K_i.
\end{array}
\end{equation*}
Here the jump term $K_3$ can be estimated as follows, recalling $g_\b(\tau,0)\equiv0$.
$$
K_3=\int_{\tau_1}^\tau\Big[\nu(\f{\t_y}{v}-\f{\T_y}{V}) Zg_\b^2\Big](\tau)d\tau
=\nu\int_{\tau_1}^\tau  g_\b^2(\tau,0)(\f{\t_y}{v}-\f{\T_y}{V})(\tau,0) \big[Z\big](\tau)d\tau\equiv 0,
$$
while the terms $K_i$ ($i=1,4,5,6$) can be directly dealt with in the same manner as
in \cite{HLM} or \cite{HWY}. To bound the term $K_2$, we make use of the mass equation $\eqref{P}_1$ to
write $K_2$ in the form
\begin{equation*}
\begin{array}{lll}
\di ~~~-(p-P)\psi_y Zg_\b^2=\f{\g P\p-(\g-1)Z}{v}Zg_\b^2\p_\tau=\di\f{\g PZg_\b^2}{2v}(\p^2)_\tau
-\f{(\g-1)Z^2 g_\b^2}{v}\p_\tau\\[4mm]
=\di\Big( \f{\g PZ \p^2g_\b^2-2(\g-1)\p Z^2g_\b^2}{2v}\Big)_\tau
-\f{\g PZ\p^2-2(\g-1)Z^2\p}{v}g_\b(g_\b)_\tau \\[3mm]
~~~\di +\f{\g PZ\p^2-2(\g-1)Z^2\p}{2v^2}g_\b^2v_\tau- \Big(\f{2(\g-1)g_\b^2\p Z}{v}+\f{\g
Pg_\b^2\p^2}{2v}\Big)Z_\tau-\f{\g g_\b^2\p^2Z}{2v}P_\tau ,
\end{array}
\end{equation*}
where all terms on the right-side hand of the above identity can be directly bounded
in the same way as in \cite{HLM} or \cite{HWY}. Therefore, we have bounded $K_2$.

Taking $\b=\f{c_0}{2}$, one can get from Lemma \ref{lemma3} that
\begin{equation}
\begin{array}{ll}
&\di \int_{\tau_1}^\tau\dashint_{\bf R} \f{e^{-\f{c_0 y^2}{1+\tau}}}{1+\tau} Z^2 dy d\tau
\leq C\d+C\|(\p,\psi,\z)(\tau_1,\cdot)\|^2+C\|(\p,\psi,\z)(\tau,\cdot)\|^2 \\[4mm]
&\quad \di +C\int_{\tau_1}^\tau\|(\p_y,\psi_y,\z_y)\dpr^2d\tau+C\int_{\tau_1}^\tau
(1+\tau)^{-\f76}\|(\p,\psi)\|^2 d\tau \\[4mm]
&\quad \di +C(\d +\eta_1)\int_{\tau_1}^\tau \dashint_{\bf
R}(1+\tau)^{-1}e^{-\f{c_0 y^2}{1+\tau}}|(\p,\z)|^2dy d\tau.
\end{array}\label{(3.49)}
\end{equation}

Now, taking $\a=c_0$ in \eqref{(3.46)} and choosing $\d $ and $\eta_1$ suitably small,
we combine \eqref{(3.46)} with \eqref{(3.49)} to obtain
the desired estimate in Lemma \ref{lemma4}. $\hfill \Box$
\vspace{2mm}

By Lemmas \ref{lemma2} and \ref{lemma4}, we conclude
$$
\begin{array}{ll}
\di N(\tau_1,\tau_2) +\int_{\tau_1}^{\tau_2}\Big\{
\|\p_y\dpr^2+\|(\psi_y,\z_y)\dpr_1^2+\|(\psi_{y\tau},\z_{y\tau})\dpr^2\Big\} d\tau\\
\qquad \di \leq CN(\tau_1)+C\int_{0}^t(1+\tau)^{-\f76}\|(\p,\psi,\z)\|^2 d\tau+C\d^\f14.
\end{array} $$
An application of Gronwall's inequality to the above inequality gives the estimate \eqref{p1} in
Proposition \ref{priori}. This completes the proof of Proposition \ref{prop1}.



\begin{thebibliography}{99}

\bibitem{BB} S. Bianchini, A. Bressan, Vanishing viscosity
solutions of nonlinear hyperbolic systems,
 Ann. of Math. (2), 161 (2005), 223-342.

\bibitem{CP} G.Q. Chen, M. Perepelitsa, Vanishing viscosity limit of the
Navier-Stokes equations to the Euler equations for compressible
fluid flow, Comm. Pure Appl. Math., 63, (2010), 1469-1504.

\bibitem{CHT} G.Q. Chen, D. Hoff,  K. Trivisa,Global solutions of the compressible Navier-Stokes
equations with large discontinuous initial data, Comm. Partial Differential Equations, 25 (2000), 2233-2257.

\bibitem{GX}
J. Goodman, Z.P. Xin, Viscous limits for piecewise smooth solutions
to systems of conservation laws, Arch. Ration. Mech. Anal., 121 (1992), 235-265.

\bibitem{H1} D. Hoff, Construction of solutions for compressible, isentropic Navier-Stokes equations
in one space dimension with nonsmooth initial data, Proc. Roy. Soc. Edinburgh (Sect. A), 103 (1986), 301-315.

\bibitem{H2} D. Hoff, Global existence for 1D, compressible, isentropic Navier-Stokes equations with
large initial data, Trans. Amer. Math. Soc., 303 (1987), 169-181.

\bibitem{H3}
 D. Hoff, Discontinuous solutions of the Navier-Stokes equations for compressible flow,
 Arch. Rational Mech. Anal., 114 (1991), 15-46.

 \bibitem{H4} D. Hoff, Global well-posedness of the Cauchy problem for the Navier-Stokes equations
 of nonisentropic flow with discontinuous initial data, J. Differential Equations, 95 (1992), 33-74.

 \bibitem{H5} D. Hoff, Discontinuous solutions of the Navier-Stokes equations for multidimensional
 flows of heat-conducting fluids, Arch. Rational Mech. Anal., 139 (1997), 303-354.

 \bibitem{H6} D. Hoff,  Global solutions of the equations of one-dimensional, compressible flow with
 large data and forces, and with differing end states, Z. Angew. Math. Phys.,49 (1998), 774-785.

\bibitem{HL}
D. Hoff, T.P. Liu, The inviscid limit for the Navier-Stokes
equations of compressible, isentropic flow with shock data, Indiana
Univ. Math. J., 38 (1989), 861-915.

\bibitem{HLM}
F.M. Huang, J. Li and A. Matsumura, Asymptotic stability of
combination of viscous contact wave with rarefaction waves for one-dimentional
compressible Navier-Stokes system, Arch. Rat. Mech. Anal., 197 (2010), 89-116.

\bibitem{HMX} F.M. Huang, A. Matsumura, Z.P. Xin, Stability
of contact discontinuities for the 1-D compressible Navier-Stokes
equations.  Arch. Ration. Mech. Anal., 179 (2006), 55-77.

\bibitem{HWY}
F.M. Huang, Y. Wang and T. Yang, Fluid Dynamic Limit to the Riemann
Solutions of Euler Equations: I. Superposition of rarefaction waves
and contact discontinuity, Kinetic and Related Models, 3 (2010), 685-728.

\bibitem{HWY1}
F.M. Huang, Y. Wang and T. Yang, Vanishing viscosity limit of the
compressible Navier-Stokes equations for  solutions to Riemann
problem, Arch. Ration. Mech. Anal., 203 (2012) 379-413.


\bibitem{HXY} F.M. Huang, Z.P. Xin, T. Yang,
Contact discontinuity with general perturbations for gas motions.
Adv. Math., 219 (2008), 1246-1297.

\bibitem{JNS}
S. Jiang, G.X. Ni and W.J. Sun, Vanishing viscosity limit to
rarefaction waves for the Navier-Stokes equations of one-dimensional
compressible heat-conducting fluids, SIAM J. Math. Anal., 38 (2006), 368-384.


\bibitem{KM} S. Kawashima, A. Matsumura, Asymptotic
stability of traveling wave solutions of systems for one-dimensional
gas motion, Comm. Math. Phys., 101 (1985), 97-127.

\bibitem{L} T.P. Liu, Shock waves for compressible Navier-Stokes
equations are stable, Comm. Pure Appl. Math., XXXIX (1986),
565-594.

\bibitem{LX} T.P. Liu, Z.P. Xin, Nonlinear stability of rarefaction
waves for compressible Navier-Stokes equations, Comm. Math. Phys.,
118 (1988), 451-465.

\bibitem{M} S.X. Ma, Zero dissipation limit to strong
contact discontinuity for the 1-D compressible Navier-Stokes
equations, J. Differential Equations,  248 (2010), 95-110.


\bibitem{MN1} A. Matsumura, K. Nishihara, On the stability
of traveling wave solutions of a one-dimensional model system for
compressible viscous gas, Japan J. Appl. Math., 2 (1985), 17-25.

\bibitem{MN2} A. Matsumura, K. Nishihara, Asymptotics toward
the rarefaction wave of the solutions of a one-dimensional model
system for compressible viscous gas, Japan J. Appl. Math., 3 (1986),
1-13.

\bibitem{MN3} A. Matsumura, K. Nishihara, Large-time behavior
 of solutions to an inflow problem in the half space
for a one-dimensional system of compressible viscous gas, Comm.
Math. Phys., 222 (2001), 449-474.

\bibitem{NYZ} K. Nishihara, T. Yang, H. Zhao, Nonlinear stability of
strong rarefaction waves for compressible Navier-Stokes equations.
SIAM J. Math. Anal., 35 (2004), 1561-1597.

\bibitem{QW1} X. Qin, Y. Wang, Stability of wave patterns to the inflow problem
of full compressible Navier-Stokes equations, SIAM J. Math. Anal. 41 (2009), 2057-2087.

\bibitem{QW2} X. Qin, Y. Wang, Large-time behavior of solutions to the inflow problem of full
compressible Navier-Stokes equations. SIAM J. Math. Anal. 43 (2011), 341-366.

\bibitem{S}
J. Smoller,  Shock Waves and Reaction-Diffusion Equations, 2nd ed., New York:
Springer-Verlag, xxii, 1994.

\bibitem{WH} H.Y. Wang, Viscous limits for piecewise smooth solutions of the p-system, J. Math. Anal. Appl.,
299 (2004), 411-432.

\bibitem{W} Y. Wang, Zero dissipation limit of the compressible heat-conducting Navier-Stokes equations
in the presence of the shock, Acta Mathematica Scientia, 28B (2008), 727-748.

\bibitem{X}
Z.P. Xin, Zero dissipation limit to rarefaction waves for the
one-dimensional Navier-Stokes equations of compressible isentropic
gases, Comm. Pure Appl. Math., 46 (1993), 621-665.

\bibitem{XZ} Z.P. Xin, H.H. Zeng, Convergence to the rarefaction
waves for the nonlinear Boltzmann equation and compressible
Navier-Stokes equations, J. Diff. Eqs., 249 (2010), 827-871.

\bibitem{Y} S.H. Yu, Zero-dissipation limit of solutions with shocks
for systems of hyperbolic conservation laws, Arch. Ration. Mech.
Anal., 146 (1999), 275-370.


\end{thebibliography}
\end{document}